\documentclass[12pt]{amsart}
\usepackage{amssymb}
\usepackage{mathtools}
\usepackage[english]{babel}
\usepackage{epsfig}
\usepackage{mathptmx}
\usepackage{amsmath,amsfonts,amssymb}
\usepackage{mathtools}
\usepackage{enumitem}
\usepackage{soul}
\usepackage{mathrsfs}
\DeclareMathAlphabet{\mathcal}{OMS}{cmsy}{m}{n} 
\usepackage{latexsym}
\usepackage{verbatim}
\usepackage{xcolor}
\usepackage[normalem]{ulem}
\usepackage{cancel}
\usepackage{tikz-cd}\usetikzlibrary{cd}
\usepackage{hyperref}
\renewcommand{\href}[2]{\relax}
\renewcommand{\url}[1]{#1}
\usepackage{marvosym}
\usepackage[text={16.5cm, 22.25cm}, asymmetric]{geometry}

\theoremstyle{plain}
\newtheorem{theorem}{Theorem}[section]
\newtheorem{corollary}[theorem]{Corollary}
\newtheorem{lemma}[theorem]{Lemma}

\newtheorem{proposition}[theorem]{Proposition}
\theoremstyle{definition}
\newtheorem{definition}[theorem]{Definition}

\numberwithin{equation}{section}

\newtheorem{example}[theorem]{Example}

\newtheorem{remark}[theorem]{Remark}

\DeclareMathOperator{\Log}{Log}

\DeclareMathOperator{\Arg}{Arg}

\DeclareMathOperator{\D}{\mathbb{D}}
\DeclareMathOperator{\RE}{Re}
\DeclareMathOperator{\IM}{Im}
\newcommand{\hide}[1]{}

\begin{document}

	\title[Divergence Rates of Orbits of Semigroups]{Divergence Rates for Orbits of Semigroups of Holomorphic Functions}
	
	\author[Maria Kourou]{Maria Kourou$^{\S}$}
	\thanks{$^{\S}$Funded by the Deutsche Forschungsgemeinschaft (DFG, German Research Foundation) – Project number: 545293643}    
	\address{Department of Mathematics, Julius-Maximilians University of Wuerzburg, 97074, Wuerzburg, Germany}
	\email{maria.kourou@uni-wuerzburg.de}   
	
	\author{Konstantinos Zarvalis}
	\address{Department of Mathematics, Aristotle University of Thessaloniki, 54124, Thessaloniki, Greece}
	\email{zarkonath@math.auth.gr}

	
	\date{}
	
	\begin{abstract}
		We study the asymptotic behavior for forward and backward orbits of one-parameter semigroups of holomorphic self-maps of the unit disk. For both types of orbits, we determine the asymptotic growth rate of hyperbolic length and compare it with the total speed of the semigroup, obtaining sharp constants that depend on the intrinsic properties  of the semigroup and on the starting point. In the process, we discover new characterizations for standard notions in the theory of semigroups. We further explore, from a potential theoretic scope, the rate by which compact sets in petals of semigroups degenerate to a boundary point in the backward dynamical setting. 
	\end{abstract}
	
	\maketitle
	
	\section{Introduction}
	One-parameter continuous semigroups of holomorphic self-maps of the unit disk $\D:=\{z\in \mathbb{C} :|z|<1\}$ have been widely studied in recent years. Their interplay with geometric function theory and operator theory, together with their connections to dynamical systems, has stimulated research interest since the original work of Berkson and Porta~\cite{berksonporta}. A plethora of recent advances in the field can be found in the monographs~\cite{abate, Booksem, shoikhetelin} and the references therein.
	
	A family $(\phi_t)_{t\geq 0}$ of holomorphic self-maps of $\D$ is called a \emph{one-parameter continuous semigroup} (or simply, a \emph{semigroup}) in $\D$ if it satisfies the following conditions:
	\begin{enumerate}
		\item $\phi_0(z)=z$ for all $z \in \D$;
		\item $\phi_{t+s}(z) = \phi_t \left( \phi_s (z) \right)$ for all $t,s \geq 0$ and $z \in \D$;
		\item $\phi_t(z) \xrightarrow{t \to 0^{+}} z$ uniformly on compact subsets of $\D$.
	\end{enumerate}
	
	Given $z\in\D$, the curve $\gamma_z:[0,+\infty)\to\D$ with $\gamma_z(t)=\phi_t(z)$ is called the \textit{forward orbit} of $z$, while $z$ is its \textit{starting point}. If for some (and equivalently, all) $t_0 >0$ the term $\phi_{t_0}$ is a conformal automorphism of $\D$, then $(\phi_t)$ is a \textit{group}. Suppose $(\phi_t)$ is not a group of (hyperbolic) rotations, an assumption that will be maintained throughout the rest of this article without always explicitly mentioning so. A remarkable feature of such semigroups is the fact that they satisfy the continuous version of the celebrated Denjoy--Wolff Theorem (see \cite[Theorem 8.3.1]{Booksem}), according to which there exists a unique $\tau\in\overline{\D}$ such that
	$$\lim\limits_{t\to+\infty}\phi_t(z)=\lim\limits_{t\to+\infty}\gamma_z(t)=\tau, \quad\text{for all }z\in\D.$$
	This unique point $\tau$ is called the \textit{Denjoy--Wolff point} of the semigroup. Whenever $\tau\in\D$, the semigroup is called \textit{elliptic} and there exists a complex number $\lambda\in\mathbb{C}$ with $\mathrm{Re}\lambda\ge0$ so that $\phi_t'(\tau)=e^{-\lambda t}$, for all $t\ge0$. On the other hand, if $\tau\in\partial\D$, $(\phi_t)$ is said to be \textit{non-elliptic} and something similar as before holds, but this time with the angular derivative. Indeed, there exists a number $\lambda\ge0$ so that $\phi_t'(\tau)\vcentcolon=\angle\lim_{z\to\tau}\phi_t'(z)=e^{-\lambda t}$, for all $t\ge0$. In particular, whenever $\lambda>0$, $(\phi_t)$ is characterized as \emph{hyperbolic}, and otherwise \emph{parabolic}. In any case, the number $\lambda$ is uniquely determined and is called the \textit{spectral value} of $(\phi_t)$.
	
	A central theme of research over the past two decades has been the boundary and asymptotic behavior of non-elliptic semigroups in $\D$, the main question being to determine the rate at which their orbits converge to the corresponding Denjoy--Wolff point. In pursuit of a unified framework, Bracci~\cite{brspeeds} introduced three quantities, the \emph{speeds of convergence} of $(\phi_t)$: the total, orthogonal, and tangential speeds. 
	The definitions in~\cite{brspeeds}, together with the questions raised in the same work, sparked a series of further developments. The second-named author \cite{Zarvalis} disproved a conjecture on the upper bound for the tangential speed of parabolic semigroups, and Cordella~\cite{Cordella} obtained an asymptotic upper bound for the tangential speed. Concerning the orthogonal speed, Bracci, Cordella, and the first-named author~\cite{brcorkou} established its asymptotic monotonicity with respect to the semigroup for a variety of cases, while Betsakos and Karamanlis~\cite{BetKar2022} generalised this result, proving via harmonic measure the monotonicity of the orthogonal speed for all non-elliptic semigroups. Further monotonicity results, obtained through symmetrization techniques, appear in the recent work~\cite{BetsChristoulou}.
	Rates of convergence for particular classes of semigroups have also been studied in~\cite{finite_shift_2025} and~\cite{Zamorano-Zarvalis_2026}, while earlier related results can be found in~\cite[Chapter~4]{shoikhet}, \cite{BetsCDM_2020, ERSY_2008, JacLevR_2011}, \cite[Chapter~7]{shoikhetelin}, and the references therein. 
	
	One of the most natural ways to measure distances in the unit disk is to endow $\D$ with the hyperbolic metric; see Subsection~\ref{sub:hypgeom} for the precise definitions. This is exactly how the total speed is defined and why it captures the full asymptotic information of the semigroup. Let $(\phi_t)$ be a non-elliptic semigroup in $\D$ with Denjoy--Wolff point $\tau \in \partial \D$ and spectral value $\lambda\ge0$. The \emph{total speed} of $(\phi_t)$ is the hyperbolic distance $d_{\D}(0, \phi_t(0))$. Bracci showed that, asymptotically, the total speed does not depend on the starting point, so the choice of $0$ in place of any other $z\in\mathbb{D}$ entails no loss of generality. Moreover, by~\cite[Proposition~6.1]{brspeeds}, it holds
	\begin{equation}\label{speed}
		\lim_{t \to +\infty} \frac{d_{\D}(0, \phi_t(0))}{t} = \frac{\lambda}{2} .
	\end{equation}
	
	From a physical point of view, the quantity $d_{\D}(0,\phi_t(0))$ is the \textit{hyperbolic displacement} of a point that starts at $0$ and travels along the orbit $\gamma_0$ until the time $t\ge0$. Consequently, the ratio $d_{\D}(0,\phi_t(0))/t$, $t>0$, may be interpreted as the \textit{average velocity} of this motion, rather than its speed. The corresponding distance of the motion is exactly the hyperbolic length $\ell_{\D}(\gamma_0;[0,t])$ of the orbit between its points $\gamma_0(0)=0$ and $\gamma_0(t)$ (for the precise definition of the hyperbolic length, we refer the reader to Subsection \ref{sub:hypgeom}). Accordingly, the \textit{average speed} of the aforementioned motion is equal to the ratio $\ell_{\D}(\gamma_0;[0,t])/t$, for $t>0$. 
	Motivated by this observation, the initial purpose of this work is to examine the quantity $\ell_{\D}(\gamma_z;[0,t])$ and the associated average speed $\ell_{\D}(\gamma_z;[0,t])/t$ based on the choice of $z\in\D$ and the type of the semigroup.

	Before passing to the statement of our first results, let us mention that the class of parabolic semigroups is further divided into two categories. Namely, if $\lim_{t\to +\infty} d_{\D}(\phi_t(z), \phi_{t+s}(z)) >0 $, for all $z \in \D$ and for all $s>0$, then the semigroup is \textit{of positive hyperbolic step}. Otherwise, it is \textit{of zero hyperbolic step}. Note that hyperbolic semigroups also have positive hyperbolic step in the above sense; see \cite[Corollary 9.3.7]{Booksem}. On the other hand, elliptic semigroups trivially have zero hyperbolic step since their Denjoy--Wolff point lies inside the unit disk.
	
	\begin{theorem}\label{thm:hyp_length}
		Let $(\phi_t)$ be a non-elliptic semigroup in $\D$. The following are true:
		\begin{itemize}
			\item[\textup{(i)}] If $(\phi_t)$ is of positive hyperbolic step (either hyperbolic or parabolic), then for each $z\in\D$ there exists a constant $c:=c((\phi_t), z)$ such that  
			$$\lim_{t\to +\infty} \frac{\ell_{\D}(\gamma_z;[0,t])}{t}= c .$$ 
			\item[\textup{(ii)}]  If $(\phi_t)$ is of zero hyperbolic step (necessarily parabolic), then 
			$$\lim_{t \to +\infty} \frac{\ell_{\D}(\gamma_z; [0,t])}{t} = 0 \quad\text{and}\quad\liminf_{t\to +\infty} \frac{\ell_{\D}(\gamma_z; [0,t])}{\log t}\geq \frac{1}{4},$$
			for all $z\in\D$.
		\end{itemize}
	\end{theorem}
	The constant $c$ depending on the semigroup and the starting point is calculated explicitly in the proof of Theorem \ref{thm:hyp_length} in Section \ref{sec:Hyp_Length} and provides a firm grasp on how the rate is influenced. Theorem \ref{thm:hyp_length} further provides optimal constants for hyperbolic semigroups, depending on the starting point, for which forward orbits are \emph{quasi-geodesic} curves, an important notion in \textit{Gromov's hyperbolicity theory}, as we see in Remark \ref{rem:quasi-geodesics}. Elliptic semigroups are excluded from Theorem \ref{thm:hyp_length}, since the hyperbolic length $\ell_{\D}(\gamma_z;[0,+\infty))\vcentcolon=\lim_{t\to+\infty}\ell_{\D}(\gamma_z;[0,t])$, for any point $z\in\D$, is finite, as we prove in Proposition \ref{prop:elliptic}.
	
	In addition to the average speeds, Theorem \ref{thm:hyp_length} and Proposition \ref{prop:elliptic} directly yield a new characterization of the positivity of the hyperbolic step for semigroups.
	\begin{corollary}\label{cor:hyperbolic step}
		A semigroup $(\phi_t)$ in $\D$ is of zero hyperbolic step if and only if $\lim_{t\to+\infty}\frac{\ell_{\D}(\gamma_z;[0,t])}{t}=0$ for some (equivalently all) $z\in\D$.
	\end{corollary}
	Therefore, the ``average hyperbolic speed'' characterizes the hyperbolic step of semigroups. On the other hand, in \cite[Proposition 6.1]{brspeeds} Bracci shows that the ``average hyperbolic velocity'' (known as the \textit{divergence rate} in literature) characterizes whether a non-elliptic semigroup is hyperbolic or parabolic. As a result, these two quantities together provide a complete blueprint for the different types of semigroups.
	
	A significant counterpart to the forward theory is the \emph{backward dynamics}. Let $(\phi_t)$ be a semigroup in $\D$ and $z\in\D$. The inverse function $\phi_t^{-1}$ is defined at $z$ for every $t$ in a maximal interval of the form $[0,T_z)$, $T_z\in(0,+\infty]$. The set $\mathcal{W}\vcentcolon=\{z\in\D\vcentcolon T_z=+\infty\}$ is called the \emph{backward invariant set} of the semigroup and is exactly the subset of the unit disk on which the backward flow of the semigroup is well defined. Each open connected component of $\mathcal{W}$ is called a \textit{petal} of $(\phi_t)$ and these are precisely the sets on which the semigroup restricts to a group of automorphisms. If $T_z=+\infty$, the continuous curve $\tilde{\gamma}_z\vcentcolon[0,+\infty)\to\D$ with $\tilde{\gamma}_z(t)=\phi_t^{-1}(z)=\vcentcolon\phi_{-t}(z)$ is the \textit{backward orbit} of $z$. More information on backward dynamics follows in Subsection \ref{sub:semigroups}.
	
	A striking property of petals is that the continuous Denjoy--Wolff Theorem applies to the backward setting, as well.
	In fact, given a petal $\Delta$, there exists a unique fixed point $\sigma \in \partial \D$ of the semigroup such that every backward orbit $\tilde{\gamma}_z$ with starting point $z\in \Delta$ converges to $\sigma$ as $t \to +\infty$. 
	We call $\sigma$ the \emph{$\alpha$-point} of the petal. Petals fall into two classes depending on their $\alpha$-point: a petal is called \emph{parabolic} if its $\alpha$-point coincides with the Denjoy--Wolff point of the semigroup, and \emph{hyperbolic} otherwise. The angular derivative at $\alpha$-points $\sigma$ of hyperbolic petals $\Delta$ is given by $\phi_t^{\prime}(\sigma)=e^{-\mu t}$, where $\mu \in (-\infty,0)$ is called the \emph{spectral value} of $\sigma$. We will refer to it as the \emph{spectral value of the petal}. 
	For further recent advances on the properties and geometry of petals, we refer the interested reader to~\cite{BCDMG_petals, ESZ2008, ADP_petals, Compact_petals}. 
	
	Let $(\phi_t)$ be a semigroup in $\D$ and $\Delta$ a petal with $\alpha$-point $\sigma\in \partial \D$ and spectral value $\mu<0$. It is natural to ask at what rate the backward orbits converge to $\sigma$. Indeed, for every $z \in \Delta$, it holds
	\begin{equation}
		\lim_{t\to +\infty} \frac{d_{\Delta}(z, \phi_{-t}(z))}{t} = -\frac{\mu}{2} ;
	\end{equation}
	see~\cite[Lemma~13.5.1]{Booksem}. Rates of convergence for the backward flow of a semigroup have been studied in~\cite{speeds_petals}, where the notions of speeds of convergence are extended to the backward setting and their asymptotic behavior is analyzed. We also mention that Betsakos, Cruz-Zamorano, and the second-named author give, in~\cite{BCZZ_25}, a complete characterization of the orthogonal speed for all kinds of forward and backward orbits of a semigroup.
	
	Our results on hyperbolic length extend to the backward dynamics of semigroups as well. We study the asymptotic behavior of the hyperbolic length $\ell_{\D}(\tilde{\gamma}_z;[0,t])$ along a backward orbit $\tilde{\gamma}_z$ with starting point $z$ in some petal of the associated semigroup. We obtain results for both hyperbolic and parabolic petals.
	
	\begin{theorem}\label{thm:hyplength_backward}
		Let $(\phi_t)$ be a semigroup in $\D$ with a petal $\Delta$.
		For each $z\in\Delta$, there exists some constant $C:=C( (\phi_t), z, \Delta)$ such that  
		$$\lim_{t\to +\infty }\frac{\ell_{\D}(\tilde{\gamma}_z;[0,t])}{t}=C \,.$$
	\end{theorem}
	
	Once more, the constant is computed explicitly in the respective proof in Section \ref{sec:Hyp_Length}. We refrain frow writing it here to better demonstrate its dependence on the various parameters later on.
	
	A natural question is what happens with the backward orbits of points lying inside the backward invariant set $\mathcal{W}$ of a semigroup $(\phi_t)$, but not in any petal. If $z\in\D$ is such a point, we necessarily have $z\in\partial\mathcal{W}\cap\D$. The backward orbits of such point are characterized as \textit{non-regular}. In Propositions \ref{prop:non-reg} and \ref{prop:elliptic non-reg}, we will see that the average hyperbolic speed along non-regular backward orbits diverges to $+\infty$. Along with Theorem \ref{thm:hyplength_backward}, these results provide a new characterization of non-regular backward orbits through the following direct corollary.
	\begin{corollary}\label{cor:non-reg}
		Let $(\phi_t)$ be a semigroup in $\D$ with backward invariant set $\mathcal{W}$. Let $z\in \mathcal{W}$. Then, $\tilde{\gamma}_z$ is non-regular if and only if $\lim_{t\to+\infty}\frac{\ell_{\D}(\tilde{\gamma}_z;[0,t])}{t}=+\infty$.
	\end{corollary}
	
	Continuing our study in the backward dynamical setting of semigroups and motivated by problems arising from physics, we investigate the asymptotic behavior of a compact set moving backwards under the influence of a semigroup. Let $(\phi_t)$ be a semigroup in $\D$ containing some petal $\Delta$. Suppose $K \subset \overline{\Delta}$ is a \textit{non-polar} compact set (i.e. of positive logarithmic capacity) that also lies in the backward invariant set of $(\phi_t)$, which means that the backward orbits for all points in $K$ are well-defined. Pick $t>0$ large enough such that $K\cap\phi_{-t}(K) =\emptyset$. Then the triple $(\Delta,\{ K, \phi_{-t}(K)\}, \{1,0\})$ forms a so-called \textit{generalised condenser} with \textit{plates} $ K$ and $\phi_{-t}(K)$ and \textit{field} $\Delta \setminus (K\cup \phi_{-t}(K))$; see \cite[Chapter 1.3]{dubininbook} for the formal definitions. As $t\to +\infty$, the plate $\phi_{-t}(K)$ moves and eventually degenerates to a fixed point of $(\phi_t)$, due to the Denjoy--Wolff Theorem in the backward dynamics. A natural way to quantify the asymptotic behavior of the generalised condenser with the moving plate is to examine the rate at which its \textit{capacity} tends to zero. The capacity of the aforementioned generalised condenser is given by the Dirichlet integral  $\iint_{\Delta\setminus (K\cup \phi_{-t}(K))} | \nabla u|^2 dA $, where $u$ is the solution of the mixed Dirichlet-Neumann problem: 
	\begin{itemize}
		\item[(1)] $\Delta u=0$ on $\Delta\setminus (K \cup \phi_{-t}( K))$;
		\item[(2)]  $u=0$ on $\partial \phi_{-t}( K)$,\, $u=1$ on $\partial K$;
		\item[(3)] $\frac{\partial u}{\partial n}=0$ on $\partial \Delta$
	\end{itemize}
	(for further details, we refer to  \cite[Theorem 4-5]{ahlforsconformal}). However, the mixed Dirichlet-Neumann problem is closely tied to the notion of \textit{extremal length}. In particular, the extremal distance $\lambda_{\Delta}(K, \phi_{-t}(K))$ between the plates is equal to the reciprocal of the capacity of the generalised condenser in question. We refer to Subsection \ref{sub:extremallength} for the detailed definitions.  
	
	The assumption that $K$ is non-polar stems again from the physical interpretation. If $K$ is a polar set, then $\phi_t(K)$ is also polar and the extremal distance $\lambda_{\Delta}(K, \phi_t(K))$ is trivially equal to $+\infty$ (or equivalently, the capacity of the generalised condenser is trivially equal to $0$). Thus, without loss of generality, we suppose that $K$ is non-polar and connected.
	The capacity of (generalised) condensers with moving plates under the action of semigroups, as well as its asymptotic behavior has been examined in various settings in ~\cite{kourouspectral, Compact_petals}. 
	In the present work, we extend the results of \cite{kourouspectral} to the backward dynamical setting of a semigroup. 
	
	\begin{theorem}\label{thm:extremal_length}
		Let $(\phi_t)$ be a semigroup in $\D$. Suppose $\Delta$ is a petal of $(\phi_t)$ and denote by $\mathcal{W}$ the backward invariant set of $(\phi_t)$. Let $K\subseteq\overline{\Delta} \cap \mathcal{W}$ be a compact non-polar set. The following hold:
		\begin{itemize}
			\item[\textup{(i)}] If $\Delta$ is hyperbolic with spectral value $\mu \in(-\infty,0)$, then 
			$$ \lim_{t\to +\infty} \frac{\lambda_{\Delta} (K, \phi_{-t}(K))}{t} =-\frac{\mu}{\pi}.$$
			\item[\textup{(ii)}] If $\Delta$ is parabolic, then 
			$$\lim_{t\to +\infty} \frac{\lambda_{\Delta}(K, \phi_{-t}(K))}{\log t} =\frac{2}{\pi}.$$
		\end{itemize}       	
	\end{theorem}
	
	As stated in Theorem \ref{thm:hyplength_backward} and in detail in Subsection \ref{subsec:hypback}, the asymptotic behavior of the hyperbolic length depends on the type of the semigroup that contains the parabolic petal. In contrast, Theorem \ref{thm:extremal_length} shows that, in the case of parabolic petals, the asymptotic behavior of the extremal length is independent of both the hyperbolic step of the semigroup and the location of the compact set $K$. This highlights, in fact, a significant difference between hyperbolic and extremal length. 
	
	However, the situation gets complicated when $K$ does not lie in the closure of some petal, but on the boundary of the backward invariant set. In this setting, the extremal distance $\lambda_{\D}(K, \phi_{-t}(K))$ gets infinite faster than in the cases covered by Theorem \ref{thm:extremal_length}, as $t \to +\infty$, and its rate depends on other intrinsic properties of $(\phi_t)$, as we see in Section \ref{sec:non-regular}.  
	
	The paper is structured in the following way. Section \ref{background} is devoted to the presentation of all auxiliary material for our proofs, providing a deeper insight to the theory of one-parameter semigroups in $\D$ and the properties of hyperbolic geometry and extremal length. Theorems \ref{thm:hyp_length} and \ref{thm:hyplength_backward} are proved in Section \ref{sec:Hyp_Length}. Section \ref{sec:Extremal_length} contains the proof of Theorem \ref{thm:extremal_length}, whereas in Section \ref{sec:non-regular}, we examine the rate of convergence for the hyperbolic length along non-regular backward orbits of a semigroup as well as the speed of convergence for compact sets that lie on the boundary of the invariant set.

	\section{Preparation for the Proofs}\label{background}
	
	\subsection{Hyperbolic Geometry}\label{sub:hypgeom}
	We start with a brief rundown of the theory of hyperbolic geometry we are going to need for our purposes. Given a simply connected domain $\Omega\subsetneq\mathbb{C}$, we denote by $\lambda_{\Omega}(z)$ the \textit{density} at $z\in\Omega$ of the hyperbolic metric of $\Omega$. Specifically for the unit disk, we have $ \lambda_{\D} (z)= (1-|z|^2)^{-1}$, for $z \in \D$. The classical Schwarz--Pick Lemma may be rewritten in the context of hyperbolic geometry as follows. 
	\begin{lemma}[Schwarz-Pick Lemma]\cite[Theorem 6.4]{bearmin}
		Suppose $\Omega_1, \Omega_2$ are simply connected proper subregions of $\mathbb{C}$ and that $f: \Omega_1 \to \Omega_2$ is holomorphic. Then 
		$$\lambda_{\Omega_2}(f(z)) | f^{\prime}(z)| \leq \lambda_{\Omega_1}(z),$$
		with equality holding if and only if $f$ maps $\Omega_1$ conformally onto $\Omega_2$. 
	\end{lemma}
	
	Let $\gamma:(\alpha,\beta)\to\D$ be a piecewise smooth curve where $-\infty\le\alpha<\beta\le+\infty$. Let $t_1, t_2 \in (\alpha,\beta)$ with $t_1<t_2$. The \emph{hyperbolic length} $\ell_{\D}(\gamma;[t_1,t_2])$ of $\gamma$ between its points $\gamma(t_1)$ and $\gamma(t_2)$ is given by
	$$\ell_{\D}(\gamma;[t_1,t_2])=\int\limits_{t_1}^{t_2}\lambda_{\D}(\gamma(t))|\gamma'(t)|dt.$$
	In a similar manner, the \textit{hyperbolic distance} between $z,w \in \D$ is given through
	$$d_{\D}(z,w)=\inf\limits_{\gamma}\int\limits_{\gamma}\lambda_{\D}(\zeta)|d\zeta|,$$
	where the infimum is taken over all piecewise smooth curves $\gamma$ joining $z$ to $w$ inside the unit disk. 
	Given such a piecewise smooth curve $\gamma:[0,1]\to\D$ with $\gamma(0)=z$ and $\gamma(1)=w$, it is true that $d_{\D}(z,w)\le \ell_{\D}(\gamma;[0,1])$, with equality if and only if $\gamma$ is a \textit{hyperbolic geodesic}. The hyperbolic geodesics in $\D$ are the diameters of $\D$ and all circular arcs in $\D$ that are perpendicular to the unit circle. 
	In fact, the hyperbolic distance can be calculated as follows:
	$$d_{\D}(z,w)=\frac{1}{2}\log\frac{|1-\overline{z}w|+| z-w|}{|1-\overline{z}w|-| z-w|} \, ; \quad z, w\in \D. $$
	
	In what follows, we cite auxiliary lemmas providing an explicit formula for the hyperbolic density in certain simply connected domains. 
	
	\begin{lemma}\label{lem:density_strip}
		Let $S_x:=\{w \in \mathbb{C} : x<\textup{Im} w< x+d\}$ be a horizontal strip of width $d>0$, where $x \in \mathbb{R}$. Then the density of the hyperbolic metric is given by 
		$$\lambda_{S_x}(z)=\frac{\pi}{2d} \frac{1}{\cos\left[\frac{\pi}{d} (\textup{Im}z-x-\frac{d}{2} ) \right]} \, , \quad z\in S_x  .$$ 
		\begin{proof}
			Let $x\in \mathbb{R}$ and $d>0$. Adjusting \cite[Example 7.9]{bearmin} and using the function $z\mapsto \frac{\pi}{d}(z-ix) -\frac{\pi}{2}i$, which maps $S_x$ conformally onto the horizontal strip $\{ z\in \mathbb{C} : |\textup{Im} z| <\frac{\pi}{2}\} $, we obtain the desired formula for $\lambda_{S_x}$.     
		\end{proof}
	\end{lemma}
	
	\begin{lemma}[{cf. \cite[Example 7.2]{bearmin}}]\label{lem:density_halfplane}
		Let $H=\{w \in \mathbb{C}: \textup{Im} w> \alpha\}$, $\alpha \in \mathbb{R}$, be a horizontal half-plane. Then, for the density of the hyperbolic metric we have $$\lambda_{H}(z)=\frac{1}{2(\textup{Im} z -\alpha)} \, , \quad z\in H. $$
	\end{lemma}
	\begin{remark} \label{rem:density_halfplane}
		In case $H=\{w \in \mathbb{C}: \textup{Im} w< \alpha\}$, $\alpha \in \mathbb{R}$, it follows $$\lambda_{H}(z)=\frac{1}{2(\alpha - \textup{Im} z )} \, , \quad z\in H. $$ So, in general, $\lambda_{H}(z)=\frac{1}{2|\textup{Im} z -\alpha|}. $
	\end{remark}
	\begin{lemma}[{cf. \cite[Example 7.7]{bearmin}}]\label{lem:density_koebe}
		Let $K:=\mathbb{C} \setminus(-\infty,0]$ be a Koebe-like domain. The hyperbolic density is given by $$\lambda_K(z) = \frac{1}{4|z| \cos \left(\frac{\textup{arg} z}{2} \right)}\, ,  \quad z \in K \,. $$
	\end{lemma}
	Inside a simply connected domain, we can relate the hyperbolic density to the distance from the boundary. Let $\Omega$ be a simply connected, proper subdomain of $\mathbb{C}$ and define $\delta_\Omega(z) := \text{dist}(z, \partial \Omega) $ to be the Euclidean distance of a point $z\in \Omega$ from the boundary of $\Omega$. Then, for any point $z\in \Omega$,
	\begin{equation}\label{eq:relationquasi}
		\frac{1}{4\delta_\Omega (z)} \leq \lambda_\Omega (z) \leq \frac{1}{\delta_\Omega (z)}  \, .
	\end{equation}

	\subsection{ One-Parameter Semigroups in $\D$}\label{sub:semigroups}
	
	As observed also in the Introduction, the iterates of one-parameter semigroups in $\D$ share similar properties. In particular, this is true for fixed points. 
	
	\begin{definition}\cite[Chapter 12]{Booksem}
		\normalfont
		Let $(\phi_t)$ be a semigroup in $\D$. A point $\sigma\in\partial\mathbb{D}$ is called a \textit{boundary fixed point} of $(\phi_t)$ provided $\angle\lim\limits_{z\to\sigma}\phi_t(z)=\sigma$, for all $t\ge0$. A boundary fixed point of $(\phi_t)$ is characterized as \textit{regular}, if the angular derivative of $\phi_t$ at $\sigma$ is finite, for all $t \geq 0$. Otherwise, we say that it is \textit{non-regular}. 
		
		Any regular boundary fixed point of a semigroup, other than the Denjoy--Wolff point, is called \textit{repelling}. If $\sigma$ is a repelling fixed point of $(\phi_t)$, then there exists a $\mu\in(-\infty,0)$ such that the angular derivative $\phi_t^{\prime}(\sigma)\vcentcolon=\angle\lim_{z\to\sigma}\phi_t'(z)=e^{-\mu t}$, for all $t\ge0$. This negative real number $\mu$ is called the \textit{repelling spectral value} of $(\phi_t)$ at $\sigma$. 
		Any non-regular boundary fixed point is characterized as \textit{super-repelling}. 
	\end{definition}
	
	If $(\phi_t)$ is a group in $\D$, then every iterate of $(\phi_t)$ is a conformal automorphism of $\D$. Groups of automorphisms of $\D$ with a fixed point inside the unit disk are called \emph{elliptic groups} and as stated in the Introduction, they are excluded from the scope of this study.
	
	For every non-elliptic semigroup $(\phi_t)$ in $\D$, there exists a unique Riemann mapping $h$ (up to some normalization; usually we consider $h(0)=0$ and $h^{\prime}(0)>0$) such that 
	$$h(\phi_t(z))=h(z)+t, \quad \text{for all }z\in\mathbb{D} \text{ and }t\ge0.$$
	The mapping $h$ is called the \emph{Koenigs function} of $(\phi_t)$ and the simply connected domain $\Omega:=h(\mathbb{D})$ is called the \textit{Koenigs domain} of the semigroup. The Koenigs domain $\Omega$ of a non-elliptic semigroup is \textit{convex in the positive direction}, which means that $\{ w+t: t\ge0 \} \subset \Omega$, for every $w \in \Omega$. By its definition, it can be understood that a major benefit of the Koenigs function is the linearization of the orbits of the points in $\D$ under the semigroup.  
	Basically, the image through $h$ of the forward orbit of some $z\in\mathbb{D}$ is a horizontal half-line emanating from $h(z)$ towards $\infty$ in the positive direction (i.e. with constant imaginary part and increasing real part). A remarkable attribute is that the reverse is also true. For any convex domain $\Omega$ in the positive direction, we may take a Riemann mapping $g_{\Omega}:\D\to \Omega$ and construct a non-elliptic semigroup $(\psi_t)$ in $\D$ via the relation $\psi_t(z)\vcentcolon = g_{\Omega}^{-1}(g_{\Omega}(z)+t)$, for any $t\ge 0$ and $z\in \D$. 
	
	The shape of the Koenigs domain leads to a classification within the class of semigroups in $\D$. As a matter of fact, a semigroup is hyperbolic if and only if its Koenigs domain $\Omega$ is contained inside a horizontal strip of width $\frac{\pi}{\lambda}$, where $\lambda$ is the spectral value of the semigroup. Furthermore, the semigroup is parabolic of positive hyperbolic step if and only if $\Omega$ is contained inside a horizontal half-plane, but not in a horizontal strip. In any other case, the semigroup is parabolic of zero hyperbolic step.
	
	The Koenigs function can also be uniquely defined for elliptic semigroups that are not groups. In this case, $h$ satisfies
	\begin{equation}\label{mapofhelliptic}
		h(\phi_t(z)) = e^{-\lambda t} h(z),\quad \text{for all }  z\in \D \text{ and } t\ge0, 
	\end{equation}
	where $\lambda\in\mathbb{C}$ is the spectral value of $(\phi_t)$, and $h(\tau)=0$, where $\tau\in\D$ is the Denjoy--Wolff point of $(\phi_t)$. The Koenigs domain $\Omega$ is \textit{$\lambda$-spirallike} with respect to $0$, since $0\in\Omega$ and $e^{-\lambda t }\Omega \subseteq \Omega$, for all $t \geq 0$. 
	
	The Koenigs function associated to an elliptic semigroup maps the forward orbit of $z\in\mathbb{D}$ onto a half-spiral, 
	which winds around $h(\tau)=0$ infinitely many times. Further information on the Koenigs function of a one-parameter semigroup and advances on its geometry and overall behavior can be found in \cite[Chapter 9]{Booksem} as well as in references therein. 
	
	For our proofs, we require some further information about the forward orbits of a semigroup. In particular, when $(\phi_t)$ is a hyperbolic semigroup with Denjoy--Wolff point $\tau\in \partial \D$, all forward orbits converge to $\tau$ non-tangentially. In fact, for every point $z\in \D$, the slope by which its forward orbit $\gamma_z$ converges to $\tau$ is defined by:
	$$\theta(z):= \textup{Slope}[\gamma_z(t),\tau] = \lim_{t\to +\infty} \textup{Arg}(1-\overline{\tau} \phi_t(z)) \in \left(-\frac{\pi}{2}, \frac{\pi}{2}\right) \, , $$
	where $\textup{Arg}$ denotes the principal argument. In addition, according to \cite[Theorem 17.3.3]{Booksem}, every forward orbit $\gamma_z$, $z\in \D$, is a \textit{quasi-geodesic} in $\D$, meaning that there exist constants $A\geq 1$ and $B\geq0$ such that for every $0\leq s\leq r<+\infty$
	$$ \ell_{\D}(\gamma_z;[s,r]) \leq A d_{\D}(\gamma_z(s), \gamma_z(r)) +B.$$
	For parabolic semigroups of positive hyperbolic step, we have that either $\theta(z)=\frac{\pi}{2}$, for all $z\in \D$, or $\theta(z)=-\frac{\pi}{2}$, for all $z\in \D$, and thus, all forward orbits converge tangentially to the Denjoy--Wolff point. For parabolic semigroups of zero hyperbolic step, the situation is not clear, as we can even have cases where the slope $\theta(z)$ is an interval, for some $z\in \D$ and thus, forward orbits can oscillate (see e.g. \cite[Proposition 17.6.1]{Booksem}). 
	
	Passing to the backward dynamical setting of a semigroup $(\phi_t)$ in $\D$, the \emph{backward invariant set} is given by $\mathcal{W} : =\bigcap_{t\geq 0} \phi_t(\D)$. Trivially, if $(\phi_t)$ is a non-elliptic group, $\mathcal{W}=\D$. Therefore, when it comes to backward dynamics, we may exclude groups, due to their trivial behavior. Every non-empty connected component $\Delta$ of the interior of $\mathcal{W}$ is called a \textit{petal} (see \cite[Definition 13.4.1]{Booksem}). Petals are non-empty simply connected subdomains of $\D$. The restriction of every $\phi_{-t}$ in a petal $\Delta$ is an automorphism (i.e. $\phi_t(\Delta)=\Delta)$. Moreover, $\tau \in \partial \Delta$ and for every $z \in \Delta $, its backward orbit $\tilde{\gamma}_z$ is \emph{regular} as we mentioned in the Introduction. The regularity of $\tilde{\gamma}_z$ is also equivalent to $\sup \{d_{\D}(\tilde{\gamma}_z(t), \tilde{\gamma}_z(t+1)): t\geq 0\} < +\infty$. Without loss of generality, by denoting $\phi_{-t}(z):= \phi_t^{-1}(z)$, for all $z\in \Delta$, we can extend the semigroup in such a way that the parameter becomes real-valued. 
	This extension is further applied to the images through the Koenigs function $h$. It can be easily seen that given $z\in\mathcal{W}$, in the case of non-elliptic semigroups, the image through $h$ of the backward orbit with starting point $z$ is actually a horizontal half-line emanating from $h(z)$ towards $\infty$ in the negative direction. In the case of elliptic semigroups, such a backward orbit is mapped through $h$ to a half-spiral reaching towards $\infty$.
	
	All regular backward orbits of a semigroup converge to a regular boundary fixed point, which is either the Denjoy--Wolff point or a repelling fixed point of the semigroup. Furthermore, all regular backward orbits of the same petal converge to the same regular fixed point. The petal is called \emph{parabolic}, if the point of convergence is $\tau$. Parabolic petals can only occur when the semigroup is itself parabolic. Otherwise, the petal is called \emph{hyperbolic} and the point to which all regular orbits converge in this petal is a repelling fixed point of the semigroup. This point is said to be the \textit{$\alpha$-point} of the petal. In fact, there is a one-to-one correspondence between repelling fixed points and hyperbolic petals. 
	
	Let $(\phi_t)$ be a non-elliptic semigroup in $\D$ with Koenigs function $h$ and Koenigs domain $\Omega$. The image of a parabolic petal $\Delta$ is a \textit{maximal} horizontal half-plane in $\Omega$, in the sense that there exists no other horizontal half-plane $H$ such that $h(\Delta) \subset H \subset \Omega$. When it comes to a hyperbolic petal $\Delta$, its image under the associated Koenigs function $h$ is a \textit{maximal} horizontal strip of width $-\frac{\pi}{\mu}$, where $\mu$ is the repelling spectral value of its $\alpha$-point. In the case where $(\phi_t)$ is an elliptic semigroup with Denjoy--Wolff point $\tau \in \D$ and spectral value $\lambda$, the image $h(\Delta)$ of a petal $\Delta$ is a maximal $\lambda$-spirallike sector in $\Omega$ of center $e^{i \theta_0}$, for some $\theta_0 \in [-\pi,\pi)$, and amplitude $2a:=-\frac{|\lambda|^2 \pi}{\mu \RE \lambda}$, where $\mu$ is the repelling spectral value of the $\alpha$-point of $\Delta$. More concretely, 
	$$h(\Delta) = \textup{Spir}\left[\lambda , 2a, \theta_0 \right]: = \{e^{t \lambda+i \theta} : t \in \mathbb{R}, \theta \in (-a+\theta_0, a+\theta_0)\}. $$
	Every petal of an elliptic semigroup is characterized as \textit{hyperbolic}. The similarity between the two types of hyperbolic petals we described relies on the fact that both of them contain a repelling fixed point in their boundaries.
	
	In our proofs, the following lemma is required to handle the case of a hyperbolic petal of an elliptic semigroup similarly to those of non-elliptic semigroups. The proof comprises of simple calculations, so we omit it for the sake of brevity.
	\begin{lemma}\label{lem:mapofpetal}\cite[Lemma 2.2]{speeds_petals}
		Suppose that $(\phi_t)$ is an elliptic semigroup in $\D$ with spectral value $\lambda \in \mathbb{C} $ and Koenigs function $h$. Let $\Delta$ be a petal of $(\phi_t)$ corresponding to the repelling fixed point $\sigma$ with repelling spectral value $\mu$. Consider the function $g(z):= \frac{e^{-i \Arg \lambda }}{\cos(\Arg \lambda)} \Log z+ i (a -\theta_0)$, where $a:=-\frac{|\lambda|^2 \pi}{2 \mu \RE \lambda}$, $\theta_0 \in [-\pi,\pi)$, and $e^{i \theta_0}$ is the center of $h(\Delta)$. Then $g$ maps $h(\Delta)$ onto the horizontal strip $\{z\in \mathbb{C}: 0<\IM z<2 a\}$. 
	\end{lemma}
	
	Last but not least, we state some useful information on non-regular backward orbits. A non-regular backward orbit either lies on the boundary of a petal, in which case it converges tangentially to the corresponding $\alpha$-point, or it converges to a super-repelling fixed point. In the latter case, the convergence can be tangential, non-tangential or oscillating. Let $\tilde{\gamma}_z:[0,+\infty)\to\mathbb{D}$ be a non-regular backward orbit for a semigroup $(\phi_t)$ with starting point $z\in \D$. Through the Koenigs function $h$, we move to the associated planar domain $\Omega$. Then, each such non-regular backward orbit of a semigroup in $\D$ satisfies the following behavior: 
	\begin{equation}\label{distancedegen}
		\lim\limits_{t\to+\infty}\delta_\Omega(h(\tilde{\gamma}_z(t)))=0,
	\end{equation}
	where $\delta_\Omega$ denotes the Euclidean distance from the boundary of $\Omega$. For further information on forward and backward orbits of one-parameter semigroups in $\D$, we refer the interested reader to \cite{Booksem}.

	\subsection{ Extremal Length}\label{sub:extremallength}
	A conformally invariant quantity that has a strong connection with the boundary behavior of holomorphic maps is the extremal length. Let $\Gamma$ be a family of rectifiable arcs in a domain $\Omega$ or, more generally, a finite union of such arcs. Consider all the non-negative Borel measurable functions $\rho$ on $\Omega$ such that the area
	$$\textup{A}( \rho, \Omega)\vcentcolon = \iint_{\Omega} \rho(z)^2 dxdy \in (0, + \infty).$$
	The \textit{extremal length} of $\Gamma $ is defined as 
	$$\lambda_{\Omega} ( \Gamma )\vcentcolon = \sup_{\rho} \frac{\left(\inf_{ \gamma \in \Gamma } \int_{\gamma} \rho(z) |dz| \right)^2}{\textup{A}( \rho, \Omega)}, $$
	where the supremum is taken over all the functions $\rho$ with $0< \textup{A} (\rho, \Omega) < + \infty$.
	It follows from the definition that if $\Gamma_1 $ is a family of curves in $\Omega$ contained in $\Gamma$, then $\lambda_{\Omega}(\Gamma_1)\ge \lambda_{\Omega}(\Gamma)$. 
	
	The conformal invariance of extremal length is also a direct aspect of the definition, as we can observe in \cite[p.130]{margarnett}.
	The most important instance of extremal length is the \textit{extremal distance} between sets lying in the closure of a domain. Suppose $E,F \subset \overline{\Omega}$. Then the extremal distance between $E$ and $F$ in $\Omega$ is defined as $\lambda_{\Omega} (E,F) : = \lambda_{\Omega \setminus (E \cup F)} (E,F) = \lambda_{\Omega} (\Gamma)$,
	where $\Gamma $ is the family of arcs or closed curves in $\Omega$ that join $E$ and $F$. 
	The \textit{conjugate extremal distance} between $E$ and $F$ is defined as $\lambda_{\Omega}^{\star}(E,F):= \lambda_{\Omega} ( \Gamma^{\star} ),$
	where $ \Gamma^{\star} $ is the family of all curves that separate $E$ from $F$ in $\Omega$.
	As we can see in \cite[p.491]{margarnett}, the conjugate and extremal distances satisfy the relation
	\begin{equation}\label{eq:conjugate}
		\lambda_{\Omega}(E,F) =\frac{ 1}{ \lambda_{\Omega}^{\star}(E,F)}.
	\end{equation}
	We will need the following two properties of extremal length. 
	\begin{lemma}[The extension rule]\cite[Theorem 4-1]{ahlforsconformal} \cite[p.134]{margarnett}\label{lem:extension}
		Let $\Omega$ and $\Omega^{\prime}$ be domains in $\mathbb{C}$ with $\Omega \subseteq \Omega^{\prime}$.  Let $\Gamma$ be family of curves in $\Omega $. Then $\lambda_{\Omega}(\Gamma) =\lambda_{\Omega^{\prime}}(\Gamma)$. Moreover, let $\Gamma^{\prime} $ be a family of curves in $\Omega^{\prime}$ such that every $\gamma^{\prime} \in  \Gamma^{\prime} $ contains some $\gamma \in  \Gamma $. Then $\lambda_{\Omega} ( \Gamma^{\prime} ) \geq \lambda_{\Omega} ( \Gamma )$.
	\end{lemma}
	
	\begin{lemma}\cite[Theorem 2.8]{ohtsuka}\label{non-disjointparallel}
		Let $ \Gamma_1 ,..., \Gamma_n $ be families of rectifiable curves in a domain $\Omega \subset \mathbb{C}$. Then 
		$$\lambda_{\Omega} \left( \bigcup_{k=1}^n \Gamma_k  \right)^{-1} \leq \sum_{k=1}^{n} \frac{1}{\lambda_{\Omega} (\Gamma_k )}. $$ 
	\end{lemma}
	
	Two basic examples of calculation of the extremal length are the cases of a rectangle and an annulus; see \cite[p.131-3]{margarnett}. Denote by $R$ the rectangle $\{ x+iy \in \mathbb{C} : 0<x<\ell \, , \,  0<y<h\}$ of length $\ell$ and height $h$. Then, the extremal length between its vertical sides $E$ and $F$ is 
	\begin{equation}\label{eq:rectangle}
		\lambda_{R}(E,F) =\frac{\ell}{h}.
	\end{equation}
	Suppose now that $\Omega = \{z\in \mathbb{C} : r<|z|<R\}$ for some $0<r<R$. Then $\lambda_{\Omega}(r \mathbb{T},R\mathbb{T} )= \frac{1}{2\pi} \log \frac{R}{r}$, where $\mathbb{T}$ denotes the unit circle. 
	In our setting, we will also need the extremal distance in a domain that is a half-annulus. 
	\begin{lemma}\label{lem:halfannulus}
		Let $\mathbb{H}:=\{ z\in \mathbb{C} : \textup{Im} z >0\}$ and $R>r>0$. Then $\lambda_{\mathbb{H}}(r \mathbb{T},R\mathbb{T})= \frac{1}{\pi} \log \frac{R}{r}$. 
		\begin{proof}
			Consider the conformal mapping $g:\mathbb{H} \to \mathbb{C}$ with $g(z)=\Log z$ which maps the upper half-plane conformally onto the horizontal strip $S\vcentcolon=\{z\in \mathbb{C} : \textup{Im} z\in (0, \pi)\}$. In particular, $r\mathbb{T}$ and $R\mathbb{T}$ are mapped onto the line segments $I_r:=\log r +i (0,\pi)$ and $I_R:=\log R+ i (0,\pi)$, respectively. 
			By the conformal invariance of extremal length, the conjugate extremal distance in $\mathbb{H}$ is given by $\lambda_{\mathbb{H}}^{\star}(r \mathbb{T},R\mathbb{T}) = \lambda^{\star}_S(I_r, I_R)$. The family of curves separating $I_r$ from $I_R$ in $S$ coincides with the family of curves separating $I_r$ from $I_R$ in the rectangle $(\log r,\log R)\times(0,\pi)$. As a result, \eqref{eq:rectangle} yields $\lambda_{\mathbb{H}}^{\star}(r \mathbb{T},R\mathbb{T}) = \frac{\pi}{\log R -\log r}$. Relation \eqref{eq:conjugate} provides at once the desired equality.
		\end{proof}
	\end{lemma}
	For further information on extremal length, the reader may refer to \cite[Chapter 4]{ahlforsconformal}, \cite[Chapter IV]{margarnett} and \cite{ohtsuka}.

	\section{Hyperbolic Length - Proof of Theorems \ref{thm:hyp_length} and \ref{thm:hyplength_backward}}\label{sec:Hyp_Length}
	
	We start the main part of this work with the examination on the hyperbolic length of the orbits. First, we deal with the forward setting, and then turn our attention to backward dynamics.
	
	\subsection{Forward Dynamics } 
	We start with the statements and proofs of all our results concerning the asymptotic behavior of the hyperbolic length along forward orbits.
	
	Let us first explain the case of elliptic semigroups and why the hyperbolic length along forward orbits is bounded. 
	When we have an elliptic semigroup, which is not a group, all forward orbits converge to an interior point of the unit disk and in fact, they are spiralling around it. However, in \cite{betszar_rect} the authors prove that all such forward orbits are rectifiable. In the following proposition, we prove that rectifiability of forward orbits of elliptic semigroups is preserved in the hyperbolic geometric sense, as well. 
	
	\begin{proposition}\label{prop:elliptic}
		Let $(\phi_t)$ be an elliptic semigroup in $\D$. Fix $z\in \D$. Then the hyperbolic length $\ell_{\D}(\gamma_z;[0,+\infty))$ is finite. 
		\begin{proof}
			Let $\tau\in \D$ be the Denjoy--Wolff point of $(\phi_t)$. The orbit $\gamma_z$ winds about $\tau$ and the set $\gamma_z([0, +\infty)) \cup\{\tau\}$ is a compact subset of $\D$. Hence there exists $M:= \max\{|\gamma_z(t)|: t\geq 0\}$ with $M\in (0, 1)$. For the hyperbolic length, we have 
			\begin{eqnarray*}
				\ell_{\D}(\gamma_z;[0,t]) &=& \int_0^t \lambda_{\D}(\gamma_z(s)) |\gamma^{\prime}_z(s)| ds = \int_0^t \frac{ |\gamma^{\prime}_z(s)|}{1-|\gamma_z(s)|^2} ds
				\\
				&\leq&\frac{1}{1-M^2}   \int_0^t |\gamma^{\prime}_z(s)| ds = \frac{\ell_{e}(\gamma_z;[0,t])}{1-M^2},
			\end{eqnarray*}
			where $\ell_{e}$ denotes the Euclidean length. However, \cite[Theorem 1.1]{betszar_rect} certifies the existence of a constant $L$ such that $\ell_{e}(\gamma_z;[0,t]) \leq L$, for all $t\ge0$. Since $\ell_{\D}(\gamma_z;[0,t])$ is increasing in $t\geq 0$, the limit as $t \to +\infty$ exists in $(0,+\infty]$, and in particular $\ell_{\D}(\gamma_z;[0,+\infty])) \leq \frac{L}{1-M^2}$.
		\end{proof}
	\end{proposition}

	We now turn our focus to non-elliptic semigroups. This time, forward orbits converge to a boundary point and so $\ell_{\D}(\gamma_z;[0,+\infty))$ is trivially infinite. For the convenience of the reader, we restate Theorem \ref{thm:hyp_length} with an explicit description of the constants mentioned.  
	\begin{theorem}\label{thm:hyp_length_const}
		Let $(\phi_t)$ be a non-elliptic semigroup in $\D$ with spectral value $\lambda \in [0,+\infty)$, Koenigs function $h$ and Koenigs domain $\Omega$. The following are true:
		\begin{itemize}
			\item[\textup{(i)}] If $(\phi_t)$ is a hyperbolic semigroup, then  $$\lim_{t\to +\infty} \frac{\ell_{\D}(\gamma_z;[0,t])}{t}= \frac{\lambda}{2\sin \left( \lambda \textup{Im} h(z) -\lambda x \right)} ;$$
			for any $z \in \D$, where $x=\inf \{ \textup{Im} h(\zeta) : \zeta \in \D \}$. 
			\item[\textup{(ii)}] If $(\phi_t)$ is a parabolic semigroup of positive hyperbolic step, then $$\lim_{t \to +\infty} \frac{\ell_{\D}(\gamma_z; [0,t])}{t}=\frac{1}{2|\textup{Im} h(z) -a|},$$
			for any $z \in \D$, where $H=\{w \in \mathbb{C} :  \textup{Im} w> a\}$ or $\{w \in \mathbb{C} :  \textup{Im} w< a\}$ is the smallest horizontal half-plane containing $\Omega$. 
			\item[\textup{(iii)}]  If $(\phi_t)$ is a parabolic semigroup of zero hyperbolic step, then
			$$\lim_{t \to +\infty} \frac{\ell_{\D}(\gamma_z; [0,t])}{t} = 0 \quad\text{and}\quad\liminf_{t\to +\infty} \frac{\ell_{\D}(\gamma_z; [0,t])}{\log t}\geq \frac{1}{4},$$
			for any $z\in\D$.
		\end{itemize}
		\begin{proof}
			Fix $t \geq  0$ and $z\in \mathbb{D}$. By the conformal invariance of $h$, we obtain $\ell_{\D}(\gamma_z;[0,t]) = \ell_{\Omega}(h \circ \gamma_z ; [0,t])$. It holds $h( \gamma_z(t)) = h(z)+t$ and thus, 
			\begin{equation*}
				\ell_{\Omega}(h \circ \gamma_z ; [0,t]) = \int_0^t \lambda_{\Omega}(h(\gamma_z(s))) \, |(h\circ \gamma_z)^{\prime}(s)| ds = \int_0^t \lambda_{\Omega}(h(z)+s) ds\, .
			\end{equation*}
			
			\textbf{(i) Hyperbolic Semigroup.} Suppose that $(\phi_t)$ is a hyperbolic semigroup in $\mathbb{D}$. Define $g(t):=\int_0^t \lambda_{\Omega}(h(z)+s) ds$. The function $g$ is differentiable in $(0, +\infty)$, with derivative $g^{\prime}(t) = \lambda_{\Omega}(h(z)+t)$. Let $S$ be the smallest horizontal strip containing $\Omega$, which we may write as $S=\{ w \in \mathbb{C} : x< \textup{Im} w <x+\frac{\pi}{\lambda}\}$ for some $x \in \mathbb{R}$, where $\lambda$ is the spectral value of $(\phi_t)$. The minimality of $S$ dictates that $x=\inf\{\mathrm{Im}h(\zeta)\vcentcolon\zeta\in\D\}$. Due to the geometry of the Koenigs domain and as a direct implication of \cite[Lemma 9.3.1]{Booksem}, $\lambda_{\Omega} (h(z)+t) \to \lambda_S(h(z))$, as $t\to +\infty$, which implies that $g^{\prime}(t) \to \lambda_S(h(z))$, as $t\to +\infty$. Hence by l'Hôpital's rule, we have $$\lim_{t \to +\infty} \frac{\ell_{\D}(\gamma_z;[0,t])}{t} = \lim_{t\to +\infty} \frac{g(t)}{t} = \lambda_S(h(z)) = \frac{\lambda}{2\sin \left( \lambda \textup{Im} h(z) -\lambda x\right)}\, ,$$
			where the final equality stems from Lemma \ref{lem:density_strip} and basic trigonometrical identities.

			\textbf{(ii) Parabolic semigroup of positive hyperbolic step} Following the same technique as in the hyperbolic case and utilising Lemma \ref{lem:density_halfplane} (or Remark \ref{rem:density_halfplane}, respectively), it turns out that
			$$\lim_{t\to +\infty } \frac{\ell_{\D}(\gamma_z;[0,t])}{t} = \lambda_H(h(z)) = \frac{1}{2|\textup{Im} h(z) -a|}\, . $$      
			
			\textbf{(iii)Parabolic semigroup of zero hyperbolic step} We initially prove the first part. Due to the properties of the Koenigs function and \eqref{eq:relationquasi}, we obtain 
			\begin{equation}\label{eq:hyplength1}
				\ell_{\D} (\gamma_z; [0,t]) =\int_0^t \lambda_{\Omega}(h(z)+s) ds  \leq \int_0^t \frac{ds}{\delta_{\Omega}(h(z)+s)}.
			\end{equation}
			Since $(\phi_t)$ is parabolic of zero hyperbolic step, the distance to the boundary satisfies $\delta_{\Omega} (h(z)+s) \to +\infty$, as $s\to +\infty$. Let $\epsilon >0$. There exists some $t_0 >0$ such that $\delta_{\Omega} (h(z)+s) >\frac{1}{\epsilon}$, for all $s \geq t_0$. Let $t>t_0$. Then, continuing in \eqref{eq:hyplength1}, 
			$$  \ell_{\D} (\gamma_z; [0,t]) \leq \int_0^{t_0} \frac{ds}{\delta_{\Omega} (h(z)+s) }  + \int_{t_0}^t \frac{ds}{\delta_{\Omega} (h(z)+s) } \leq c+ \int_{t_0}^t \epsilon ds = c+\epsilon (t-t_0) , $$
			where $c$ is a positive constant not depending on $t$.
			Therefore, 
			$$\limsup_{t\to +\infty} \frac{\ell_{\D} (\gamma_z; [0,t])}{t} \leq \lim_{t\to +\infty} \frac{c+\epsilon (t-t_0)}{t}  = \epsilon$$
			and since $\epsilon $ may be chosen arbitrarily small, it follows that the limit exists and $t^{-1} \ell_{\D} (\gamma_z; [0,t])\to0$, as $t\to +\infty$. 
			
			Now we move to the proof of the second part. Due to $\Omega$ being simply connected, convex in the positive direction, and not the whole complex plane, there exists at least one $p \in \partial \Omega$ such that $\{ p -s : s\geq 0\} \subseteq \mathbb{C} \setminus \Omega$. As a result, $\Omega -p \subseteq \mathbb{C}\setminus (-\infty, 0]=:K$. With the aid of Lemma \ref{lem:density_koebe} and the domain monotonicity property of the hyperbolic length, it follows that
			\begin{eqnarray*}
				\ell_{\D} (\gamma_z; [0,t]) 
				&=& \ell_{\Omega-p} (h\circ \gamma_z -p;[0,t]) \\
				&\geq&  \ell_{K} (h\circ \gamma_z -p;[0,t]) \\
				&=& \int_0^t \lambda_K (h(z)+s -p) ds \\
				&=& \int_0^t \frac{ds}{4 |h(z)+s-p| \cos \left(\frac{\textup{arg} (h(z)+s-p)}{2}\right)} \\
				&\geq & \frac{1}{4}  \int_0^t \frac{ds}{ |h(z)-p|+s}  \\
				&=& \frac{1}{4} \left( \log (|h(z)-p|+t) - \log|h(z)-p| \right),
			\end{eqnarray*}
			leading to $$ \liminf_{t\to +\infty} \frac{ \ell_{\D} (\gamma_z; [0,t])}{ \log t} \geq \lim_{t\to +\infty} \frac{\log (|h(z)-p|+t)}{ 4 \log t} = \frac{1}{4} . $$
		\end{proof}
	\end{theorem}
	\begin{remark}\label{rem:quasi-geodesics}
		As we see in \cite[Theorem 9.1.9]{Booksem}, for hyperbolic semigroups, it is known that $t^{-1}d_{\D}(\gamma_z(0), \gamma_z(t))\to \frac{\lambda}{2}$, as $t \to +\infty$. Combining it with Theorem \ref{thm:hyp_length}, we obtain $$\lim_{t\to +\infty} \frac{\ell_{\D}(\gamma_z;[0,t])}{d_{\D}(\gamma_z(0), \gamma_z(t))}=\frac{1}{\sin \left( \lambda \textup{Im} h(z) -\lambda x \right)}.$$
		Since all orbits of a hyperbolic semigroup converge to the Denjoy--Wolff point non-tangentially, they are \textit{quasi-geodesics} (see \cite[Theorem 17.3.3]{Booksem}); i.e. there exist $A\geq 1$ and $B\geq 0$ such that $\ell_{\D}(\gamma_z; [t_1,t_2]) \leq A d_{\D}(\gamma_z(t_1), \gamma_z(t_2)) +B$, for all $t_2\geq t_1 \geq 0$. Thus, for $t_1=0$ and $t_2$ large enough, Theorem \ref{thm:hyp_length} shows that $B=0$ and the optimal $A$ is $\frac{1}{\sin(\lambda \textup{Im}h(z) -\lambda x)}$.
	\end{remark}

	We provide an additional approach in terms of the slope of the orbit $\gamma_z$. Again, we consider $(\phi_t)$ to be a hyperbolic semigroup in $\mathbb{D} $ and fix a point $z \in \D$. Set $\theta_z \in (0, \pi)$ to be the angle by which the orbit of $z$ converges to the Denjoy--Wolff point of $(\phi_t)$. Let $S$ be the smallest horizontal strip containing $\Omega$. The Koenigs function $h$ is semi-conformal (in the sense of angle preservation) at the Denjoy--Wolff point; 
	see \cite[Proposition 17.4.5]{Booksem}. This implies that if $g:\D\to S$ is a Riemann map, then $g^{-1}(h(z)+t)$ converges to $\partial \D$, as $t \to +\infty$, still by the same angle $\theta_z$. We compute $\theta_z$ by mapping $S$ to the upper half-plane $\mathbb{H}=\{ w \in \mathbb{C} : \textup{Im} w>0 \}$. Consider the conformal mapping $f: S \to \mathbb{H}$ with $f(z):=\exp\{ \lambda z -i\lambda x\}$ and thus, $\textup{Arg} f(h(z)) = \lambda \textup{Im} h(z) -\lambda x$. In the upper half-plane $\mathbb{H}$, the angle of convergence of $f(h(z)+t)$ to $\infty$ is exactly the argument $\textup{Arg}  f(h(z)+t)=\textup{Arg} f(h(z))$ and since all mappings are conformal, the angles in $\mathbb{H}$, $S$, $\Omega$ and $\D$ are the same. Hence, $\theta_z=\lambda \textup{Im} h(z) -\lambda x$ and as a consequence, 
	$$ \frac{\ell_{\D}(\gamma_z;[0,t])}{d_{\D}(\gamma_z(0), \gamma_z(t))} \xrightarrow{t\to +\infty} \frac{1}{\sin \theta_z} . $$
	All in all, we obtain $$\frac{\ell_{\D}(\gamma_z; [0,t])}{t} \xrightarrow{t \to +\infty} \frac{\lambda}{2 \sin \theta_z} \quad \text{and} \quad \frac{\ell_{\D}(\gamma_z; [0,t])}{d_{\D}(\gamma_z(0), \gamma_z(t))}\xrightarrow{t \to +\infty}  \frac{1}{\sin \theta_z} . $$
	Therefore, the angle of convergence of the hyperbolic orbit determines and in a way, quantifies its property of being a quasi-geodesic. 
	
	When it comes to parabolic semigroups of positive hyperbolic step, the orbits land tangentially at the Denjoy--Wolff point and so, they are not quasi-geodesics. Therefore, the ratio $ \frac{\ell_{\D}(\gamma_z; [0,t])}{d_{\D}(\gamma_z(0), \gamma_z(t))}$ diverges to infinity, as $t\to +\infty$, for all $z\in \D$. Since the tangential convergence implies $\sin\theta_z=0$, the result for the hyperbolic semigroups trivially extends to parabolic semigroups of positive hyperbolic step as well.

	\begin{remark}
		According to \cite[Theorem 9.1.9]{Booksem} and \cite{brspeeds}, parabolic semigroups of positive hyperbolic step satisfy $t^{-1} d_{\D}(\gamma_z(0),\gamma_z(t)) \to 0$, as $t\to +\infty$, for all $z\in \D$. However, this does not yield information on the orbit, since it does not reveal the rate of convergence of the total speed; just that it is slower than the linear approach. By the previous theorem and combining with the preceding remark, the obtained information is improved, since it is proved that the ``divergence'' of the hyperbolic length of the orbit is linear. 
	\end{remark}

	\begin{remark}
		For parabolic semigroups of zero hyperbolic step, Theorem \ref{thm:hyp_length} dictates that the hyperbolic length diverges at least logarithmically. Of course, this rate is sharp since it is easily checked that if $\Omega = \mathbb{C} \setminus (-\infty, 0]$, the associated semigroup satisfies 
		$$\lim_{t\to +\infty} \frac{ \ell_{\D} (\gamma_z; [0,t])}{ \log t} =\frac{1}{4} .$$
	\end{remark}
	Nevertheless, the divergence rate is still not covered completely by the preceding theorem. We just know that the divergence is slower than linear and at the minimum as fast as logarithmic. A wide spectrum of rates remains in between. This gap is filled by the following proposition which shows that for every rate between logarithmic and linear, we may find a parabolic semigroup of zero hyperbolic step attaining it.
	\begin{proposition}\label{prop:all rates parabolic}
		Let $f:[0, +\infty) \to [0, +\infty)$ be an increasing, concave and smooth function satisfying $\lim_{t\to+\infty}f(t)=+\infty$ and $\lim_{t\to+\infty}f'(t)=0$. Suppose that there exists an absolute positive constant $M$ so that 
		\begin{equation}\label{eq:condition}
			\frac{1}{f'(t)}\le t+M, \quad  \text{for all } t\ge0\,.  
		\end{equation}
		Then, there exists a semigroup $(\phi_t)$ in $\D$ and $z\in \D$ such that 
		$$0<\liminf_{t \to +\infty} \frac{\ell_{\D}(\gamma_z;[0,t])}{f(t)} \le \limsup_{t \to +\infty} \frac{\ell_{\D}(\gamma_z;[0,t])}{f(t)}<+\infty. $$
		\begin{proof}
			By hypothesis $\frac{1}{f'}$ is continuous and increasing with $\frac{1}{f'(t)}\xrightarrow{t\to +\infty} +\infty$. A combination with \eqref{eq:condition} implies that there exists a simply connected domain $\Omega$ that is convex in the positive direction and is not contained in any horizontal half-plane, along with a point $w\in\Omega$ such that $\delta_{\Omega}(w+t)=\frac{1}{f'(t)}$, for all $t\ge0$. 
			
			The importance of \eqref{eq:condition} is due to the following observation. If $\Omega$ is a convex in the positive direction, simply connected, proper subdomain of the complex plane and $p\in\partial\Omega$, then for a fixed $w\in\Omega$ we have
			$$\delta_{\Omega}(w+t)\le |w+t-p|\le |w-p|+t=M+t.$$
			
			Therefore, $\Omega$ may be considered as the Koenigs domain of a parabolic semigroup $(\phi_t)$ with zero hyperbolic step. Denote by $h$ the corresponding Koenigs function and set $z=h^{-1}(w)$. Using the conformal invariance of the hyperbolic length and \eqref{eq:relationquasi}, we obtain
			\begin{eqnarray*}
				\ell_{\D}(\gamma_z;[0,t]) &=& \ell_{\Omega}(h \circ \gamma_z; [0,t]) = \int_0^t \lambda_{\Omega} (s+w) ds \\
				&\geq& \int_0^t \frac{ds}{4 \delta_{\Omega}(s+w)} =  \frac{1}{4}\int_0^t f'(s)ds\\
				&=& \frac{1}{4} (f(t)-f(0)).
			\end{eqnarray*}
			As a result, 
			\begin{equation}\label{lemma 1}
				\frac{\ell_{\D}(\gamma_z;[0,t])}{f(t)} \geq \frac{1}{4} - \frac{f(0)}{4f(t)}.     
			\end{equation}
			In similar fashion, and using the right-hand side of \eqref{eq:relationquasi}, we find
			\begin{equation}\label{lemma 2}
				\frac{\ell_{\D}(\gamma_z;[0,t])}{f(t)} \leq 1-\frac{f(0)}{f(t)}.
			\end{equation}
			Combining \eqref{lemma 1}, \eqref{lemma 2}, and taking limits as $t\to+\infty$, we arrive to the desired conclusion. The reason why this proposition fills the gap we mentioned before is because of the initial assumptions. Indeed, l'H\^{o}pital's rule and the fact that $\lim_{t\to+\infty}f'(t)=0$, means that $\lim_{t\to+\infty}\frac{f(t)}{t}=0$. In addition, algebraic considerations on \eqref{eq:condition} show that $\liminf_{t\to+\infty}(tf'(t))\ge1$. Again, l'H\^{o}pital's rule gives $\liminf_{t\to+\infty}\frac{f(t)}{\log t}\ge1$. So, this proposition covers every rate that is slower than linear, but at least as fast as logarithmic. 
		\end{proof}
	\end{proposition}
	
	Ending our discussion on the rate of the hyperbolic length for forward orbits, we remind that the combination of Proposition \ref{prop:elliptic} and Theorem \ref{thm:hyp_length} directly prove Corollary \ref{cor:hyperbolic step}. This leads to a characterization of the hyperbolic step for semigroups in the unit disk and in turn to a complete classification of semigroups in hyperbolic terms with the aid of \cite[Proposition 6.1]{brspeeds}. As a matter of fact:
	\begin{enumerate}
		\item[(1)] $(\phi_t)$ is elliptic if and only if $\lim_{t\to+\infty}d_{\D}(z,\phi_t(z))\in[0,+\infty)$ for some (equivalently all) $z\in\D$;
		\item[(2)] $(\phi_t)$ is hyperbolic if and only if $\lim_{t\to+\infty}\frac{d_{\D}(z,\phi_t(z))}{t}, \lim_{t\to+\infty}\frac{\ell_{\D}(\gamma_z;[0,t])}{t}\in(0,+\infty)$ for some (equivalently all) $z\in\D$;
		\item[(3)] $(\phi_t)$ is parabolic of positive hyperbolic step if and only if $\lim_{t\to+\infty}\frac{d_{\D}(z,\phi_t(z))}{t}=0$ and $\lim_{t\to+\infty}\frac{\ell_{\D}(\gamma_z;[0,t])}{t}\in(0,+\infty)$ for some (equivalently all) $z\in\D$;
		\item[(4)] $(\phi_t)$ is parabolic of zero hyperbolic step if and only if $\lim_{t\to+\infty}d_{\D}(z,\phi_t(z))=+\infty$ and $\lim_{t\to+\infty}\frac{d_{\D}(z,\phi_t(z))}{t}=\lim_{t\to+\infty}\frac{\ell_{\D}(\gamma_z;[0,t])}{t}=0$ for some (equivalently all) $z\in\D$.
	\end{enumerate}

	\subsection{Backward Dynamics} \label{subsec:hypback} In the current subsection, we prove Theorem \ref{thm:hyplength_backward}, which yields the asymptotic behavior of the hyperbolic length along backward orbits of a semigroup $(\phi_t)$ that lie in some petal of $(\phi_t)$.  
	First, we restate Theorem \ref{thm:hyplength_backward} with an explicit description of the constants mentioned. 
	
	\begin{theorem}\label{thm:hyplength_backward_const}
		Let $(\phi_t)$ be a semigroup in $\D$ with associated Koenigs function $h$ and a petal $\Delta$. The following are true. 
		\begin{itemize}
			\item[\textup{(i)}] If $(\phi_t)$ is elliptic (hence $\Delta$ is necessarily hyperbolic), then $$\lim_{t\to +\infty }\frac{\ell_{\D}(\tilde{\gamma}_z;[0,t])}{t}= -\frac{\mu}{2\cos\left[ \frac{\pi}{2a}(\textup{Im} f(h(z)) - a)\right]}, $$
			for any $z\in \Delta$, where $f$ is the conformal mapping from Lemma \ref{lem:mapofpetal}, $\mu\in(-\infty,0)$ is the spectral value of $\Delta$ and $a$ is an intrinsic quantity of $\Delta$.
			\item[\textup{(ii)}] If $(\phi_t)$ is non-elliptic and $\Delta$ is hyperbolic, then $$\lim_{t\to +\infty }\frac{\ell_{\D}(\tilde{\gamma}_z;[0,t])}{t}= - \frac{\mu}{2\sin \left( -\mu\textup{Im} h(z) + \mu x \right)},$$
			for any $z\in \Delta$, where $\mu\in(-\infty,0)$ is the spectral value of $\Delta$ and $x =\inf\{ \textup{Im} h(\zeta) : \zeta \in \Delta\}$. 
			\item[\textup{(iii)}] If $\Delta$ is parabolic, then  $$\lim_{t\to +\infty }\frac{\ell_{\D}(\tilde{\gamma}_z;[0,t])}{t}= \frac{1}{2|\textup{Im} h(z) -k| }, $$
			for any $z\in\Delta$, where $h(\Delta)$ is either the half-plane $\{w \in \mathbb{C} :  \textup{Im} w> k\}$ or $\{w \in \mathbb{C} :  \textup{Im} w< k\}$, $k\in\mathbb{R}$. 
		\end{itemize}
		\begin{proof}
			Fix $z\in \Delta$. Recall that $\tilde{\gamma}_z$ converges either to a repelling fixed point of $(\phi_t)$, if $\Delta$ is hyperbolic, or to the Denjoy--Wolff point, in case $\Delta$ is parabolic. Concerning the hyperbolic density, it is known from \cite[Theorem 1.3]{Compact_petals} that 
			\begin{equation}\label{eq:backward hyp length 0}
				\lambda_{\D}(\phi_{-t}(z)) |\phi^{\prime}_{-t}(z)| \to \lambda_{\Delta}(z), \text{ locally uniformly in } \Delta, \text{ as } t\to+\infty.
			\end{equation}
			For the rest of the proof, we denote by $\Omega$, as usual, the Koenigs domain of $(\phi_t)$. 
			Let us consider the various cases for semigroups and petals.     
			
			\textbf{(i)} Suppose $(\phi_t)$ is elliptic (not a group) with spectral value $\lambda \in \mathbb{C}$, $\mathrm{Re}\lambda>0$. The petal $\Delta$ is hyperbolic and its image under the associated Koenigs function $h$ of $(\phi_t)$ is a $\lambda$-spirallike sector of center $e^{i \theta_0}$ for some $\theta_0 \in (-\pi, \pi]$ and amplitude $2 a:=-\frac{|\lambda|^2 \pi}{\mu \text{Re} \lambda}$. Using Lemma \ref{lem:mapofpetal} and the conformal mapping $f(z) = \frac{e^{-i \Arg \lambda}}{\cos(\Arg \lambda)} \Log z + i (a-\theta_0)$, the spirallike sector $h(\Delta)$ is mapped onto the horizontal strip $S_a:=\{z :  0< \text{Im} z <2a\}$.
			Via its conformal invariance, the hyperbolic length of the backward orbit, for any $t\ge0$, is given by 
			$$\ell_{\D}(\tilde{\gamma}_z; [0,t]) = \ell_{\Omega}(h \circ \tilde{\gamma} _z ; [0,t]) = |\lambda h(z)|\int_0^t \lambda_{\Omega}(e^{\lambda s}h(z)) e^{\textup{Re} \lambda s} ds ,$$
			where we have used the fact that $h(\tilde{\gamma}_z(t))=e^{\lambda t}h(z)$. Following the same procedure as in the proof of Theorem \ref{thm:hyp_length}, set $$g(t):=\int_0^t \lambda_{\Omega}(e^{\lambda s}h(z)) e^{\textup{Re} \lambda s}  ds.$$ The function $g$ is differentiable in $(0, +\infty)$, with derivative $g^{\prime}(t) = \lambda_{\Omega}(e^{\lambda t}h(z)) e^{\textup{Re} \lambda t} $.     Since 
			\begin{eqnarray*}
				\lambda_{\D}(\phi_{-t}(z)) |\phi^{\prime}_{-t}(z)| &=&\lambda_{\Omega}(e^{\lambda t} h(z)) |h'(\phi_{-t}(z))||\phi'_{-t}(z)| \\
				&=& \lambda_{\Omega}(e^{\lambda t} h(z)) |e^{\lambda t}h^{\prime}(z)|
			\end{eqnarray*}
			and $\lambda_{\Delta}(z)=\lambda_{h(\Delta)}(h(z))|h'(z)|$, it follows through \eqref{eq:backward hyp length 0} that $\lambda_{\Omega}(e^{\lambda t} h(z)) e^{\textup{Re}\lambda t}\to \lambda_{h(\Delta)}(h(z))$ locally uniformly in $\Delta$, as $t\to +\infty$. This implies that $g^{\prime}(t) \to \lambda_{h(\Delta)}(h(z))$, as $t\to +\infty$ and hence, by l'Hôpital's rule, we have 
			\begin{equation}\label{eq:backward hyp length 1}
				\lim_{t \to +\infty} \frac{\ell_{\D}(\tilde{\gamma}_z;[0,t])}{t} = |\lambda h(z)| \lim_{t\to +\infty} \frac{g(t)}{t} = |\lambda h(z)| \lambda_{h(\Delta)}(h(z)) .   
			\end{equation}
			Therefore, using Lemma \ref{lem:density_strip} we obtain
			\begin{eqnarray}\label{eq:backward hyp length 2}
				\notag   |\lambda h(z)| \lambda_{h(\Delta)}(h(z)) &=&  |\lambda h(z)| \lambda_{S_a}(f(h(z))) |f^{\prime}(h(z))|\\
				\notag  &=& \frac{\pi  |\lambda h(z)|}{4a \cos(\Arg \lambda) |h(z)|} \frac{1}{\cos\left[ \frac{\pi}{2a}(\textup{Im} f(h(z)) - a)\right]} \\
				&=&\frac{-\mu}{2\cos \left[ \frac{\mu \textup{Re} \lambda}{|\lambda|^2} \theta_0- \frac{\mu \textup{Re} \lambda}{|\lambda|^2} \textup{Arg} h(z) - \mu \frac{\sin(\textup{Arg} \lambda)}{|\lambda|} \log |h(z)|\right]} \, . 
			\end{eqnarray}
			Plugging \eqref{eq:backward hyp length 2} into \eqref{eq:backward hyp length 1} yields the desired outcome.
			
			\textbf{(ii)} In the case where $(\phi_t)$ is non-elliptic and $\Delta$ is hyperbolic, the image $h(\Delta)$ is a horizontal strip of width $-\frac{\pi}{\mu}$. Set $x:=\inf\{ \textup{Im} h(\zeta) : \zeta \in \Delta\}$. In other words, $h(\Delta)=\{z\in\mathbb{C}\vcentcolon x<\mathrm{Im}z<x-\frac{\pi}{\mu}\}$. As previously
			$$\ell_{\D}(\tilde{\gamma}_z; [0,t]) = \ell_{\Omega}(h \circ \tilde{\gamma} _z ; [0,t]) 
			=\int_0^t \lambda_{\Omega}(h(z)-s) ds .$$
			
			Following the same procedure as in the proof of Theorem \ref{thm:hyp_length}, set $g(t):=\int_0^t \lambda_{\Omega}(h(z)-s) ds$. The function $g$ is differentiable in $(0, +\infty)$, with derivative $g^{\prime}(t) = \lambda_{\Omega}(h(z)-t) $. Since 
			\begin{eqnarray*}
				\lambda_{\D}(\phi_{-t}(z)) |\phi^{\prime}_{-t}(z)| &=&\lambda_{\Omega}(h(z)-t) |h'(\phi_{-t}(z))||\phi'_{-t}(z)| \\
				&=& \lambda_{\Omega}(h(z)-t) |h^{\prime}(z)|
			\end{eqnarray*}
			and $\lambda_{\Delta}(z)=\lambda_{h(\Delta)}(h(z))|h'(z)|$, \eqref{eq:backward hyp length 0} brings forth $\lambda_{\Omega}(h(z)-t) \to \lambda_{h(\Delta)}(h(z))$ locally uniformly in $\Delta$, as $t\to +\infty$. This implies that $g^{\prime}(t) \to \lambda_{h(\Delta)}(h(z))$, as $t\to +\infty$ and hence, by l'Hôpital's rule, we get
			\begin{equation}\label{eq:backward hyp length 3}
				\lim_{t \to +\infty} \frac{\ell_{\D}(\tilde{\gamma}_z;[0,t])}{t} = \lim_{t\to +\infty} \frac{g(t)}{t} =\lambda_{h(\Delta)}(h(z)) .   
			\end{equation}
			From Lemma \ref{lem:density_strip}, we have  $$\lambda_{h(\Delta)}(h(z)) = -\frac{\mu}{2} \frac{1}{\cos\left(-\mu \textup{Im} h(z) +\mu \nu-\frac{\pi}{2}  \right)}.$$
			Applying on \eqref{eq:backward hyp length 3}, we obtain the desired limit.
			
			\textbf{(iii)} In the case where $\Delta$ is parabolic, the semigroup $(\phi_t)$ is parabolic and the image $h(\Delta)$ is a half-plane of the form $\{w\in \mathbb{C} : \textup{Im} w>k\}$ or $\{w\in \mathbb{C} : \textup{Im} w<k\}$, for some $k \in \mathbb{R}$. Therefore, from Lemma \ref{lem:density_halfplane}, we obtain 
			$$\lambda_{h(\Delta)}(h(z))= \frac{1}{2|\textup{Im} h(z) -k|} . $$
			Following the steps utilised in the previous cases, we may conclude the proof. We omit the details for the sake of avoiding repetition.

		\end{proof}
	\end{theorem}
	\section{Extremal Length - Proof of Theorem \ref{thm:extremal_length}}\label{sec:Extremal_length}
	We proceed to our study concerning the way that backward orbits behave with respect to the extremal length. As mentioned in the Introduction, the corresponding investigation about forward orbits has been conducted in \cite{kourouspectral}. In this section, we prove Theorem \ref{thm:extremal_length}. The proof is based on the construction of families of crosscuts that separate a compact set from its backward orbit inside a petal, leading thus to estimates on the extremal distance. 
	Before we move on to the proof, we state the following lemma which contains a helpful estimate concerning extremal distance. 
	\begin{lemma}\label{lem:extremallengthcircles}
		Let $\mathbb{H} = \{z\in \mathbb{C} : \textup{Im} z >0\}$. Suppose $E\subset \overline{\mathbb{H}}$ is a compact, non-polar set with $0 \notin E$. Set $l_1:=\min \{ |w| :w\in E\} $ and $l_2:=\max \{ |w| :w\in E\} $. Then for all $t>2l_2$, 
		$$\frac{2}{\pi} \log t +\frac{2}{\pi} \log \frac{t+\sqrt{t^2-4l_2^2}}{2 l_2 t} \leq \lambda_{\mathbb{H}}(E, E-t) \leq  \frac{2}{\pi} \log t +\frac{2}{\pi} \log \frac{t+\sqrt{t^2-4 l_1^2}}{2l_1 t} + 2\pi \log^{-1} \frac{l_2}{l_1} . $$
		\begin{proof}
			For $0<r<R$ and $p\in \mathbb{C}$, set $A(p,r,R)$ the annulus $\{w\in \mathbb{C} : r<|w-p|<R\}$. 
			Notice that $E \subset \overline{A(0,l_1,l_2)}$.
			For $t>2l_2$, it holds $E-t \subset \overline{A(-t, l_1,l_2)}$, with $\overline{A(0,l_1,l_2)} \cap \overline{A(-t, l_1,l_2)} =\emptyset$. Set $C_j:=  \{ l_j e^{i\theta} : \theta \in [0,\pi]\}$ and $C_j^t:=  \{  -t +l_j e^{i\theta} : \theta \in [0,\pi]\}=C_j-t$, for $j=1,2$. 
			By the extension rule of extremal length (Lemma \ref{lem:extension}), it follows that 
			\begin{equation}\label{eq:paraboliclowerbound}
				\lambda_{\mathbb{H} }(E, E-t) \geq \lambda_{\mathbb{H}}(C_2^t, C_2).
			\end{equation}
			Set $\Gamma_t^{\star} $ the family of curves in $\mathbb{H}$ separating $E$ from $E-t$
			and $\Gamma_t^{1\star}$ the family of curves in $\mathbb{H}$ that separate $C_1$ from $C_1^t$. There exists a subfamily $\Gamma_{C}^{\star}$ of $\Gamma_t^{\star}$ that contains crosscuts with endpoint(s) in $\partial \mathbb{H}$ (the case of one endpoint is for when the crosscut goes to infinity on one side) that simultaneously separate $C_1$ from $C_1^t$. Thus, by the extension rule of the extremal length, 
			\begin{equation}\label{eq:crosscuts2}
				\lambda(\Gamma_{C}^{\star}) \geq \lambda (\Gamma_t^{1\star})=\frac{1}{\lambda_{\mathbb{H}}(C_1,C_1^t)}.
			\end{equation}
			
			The rest of the curves of $\Gamma_t^{\star}$ lie in one of the following categories:
			\begin{enumerate}
				\item $\Gamma_1$ : family of curves intersecting $C_1$ and $C_2$ (these are either crosscuts in $\mathbb{H}$ or curves that circle around $E$)	 
				\item $\Gamma_2$ : family of curves intersecting $C_1^t$ and $C_2^t$ (these are either crosscuts in $\mathbb{H}$ or curves that circle around $E-t$).
			\end{enumerate}
			
			Due to the construction and the fact that the mapping $z\mapsto z-t$, $t\ge0$, is a conformal automorphism of the upper half-plane thus preserving the extremal length, note that $\lambda(\Gamma_j) \geq \lambda_{\mathbb{H}}(C_1,C_2) =\lambda_{\mathbb{H}}(C_1^t, C_2^t) =\frac{1}{\pi} \log \frac{l_2}{l_1}$, for $j=1,2$. 		
			Following Lemma \ref{non-disjointparallel} and applying \eqref{eq:conjugate}, it holds
			\begin{eqnarray}\label{eq:parabolicupperbound}
				\lambda_{\mathbb{H}}(E, E-t)=\frac{1}{\lambda(\Gamma_t^{*})}&\leq& \frac{1}{\lambda(\Gamma_C^{\star})}+ \sum_{i=1}^2\frac{1}{\lambda(\Gamma_i)} \nonumber \\
				&\leq& \lambda_{\mathbb{H}} (C_1,C_1^t)+ \frac{2\pi}{\log \frac{l_2}{l_1}} .
			\end{eqnarray}
			As a result, all that remains is to calculate $\lambda_{\mathbb{H}}(r\mathbb{T}^{+}, -t+r\mathbb{T}^{+})$, where $\mathbb{T}^{+}:=\{z=re^{i\theta}: \theta \in [0,\pi]\}$, for some $r>0$ and $t>2r$. The M\"obius transformation $T_t:\overline{\mathbb{H}} \to \overline{\mathbb{H}}$ with 
			$$T_t(z)= \frac{2z+t-\sqrt{t^2-4r^2}}{2z+t+\sqrt{t^2 -4r^2}}$$
			is a conformal automorphism of $\mathbb{H}$. Through elementary calculations, we can verify that $$T_t(r \mathbb{T}^{+}) = \left\{\frac{2r}{t+\sqrt{t^2-4r^2}} e^{i\theta} : \theta \in [0, \pi] \right\}. $$
			Similarly, we can see that 
			$$T_t(-t+r \mathbb{T}^{+}) = \left\{\frac{t+\sqrt{t^2-4r^2}}{2r} e^{i\theta} : \theta \in [0, \pi] \right\}. $$
			Due to the conformal invariance of extremal length, 
			\begin{eqnarray}\label{eq:final eq}
				\notag   \lambda_{\mathbb{H}}(r\mathbb{T}^{+}, -t+r\mathbb{T}^{+}) &=&\lambda_{\mathbb{H}}\left(\frac{2r}{t+\sqrt{t^2-4r^2}} \mathbb{T}^{+},  \frac{t+\sqrt{t^2-4r^2}}{2r} \mathbb{T}^{+}\right)\\
				\notag  &=& \frac{1}{\pi} \log \left(\frac{t+\sqrt{t^2-4r^2}}{2r}\right)^2\\
				&=& \frac{2}{\pi} \log t +\frac{2}{\pi} \log \frac{1+\sqrt{1-4\frac{r^2}{t^2}}}{2r}. 
			\end{eqnarray}
			Returning back to $\lambda_{\mathbb{H}}(E, E-t)$ in relations \eqref{eq:paraboliclowerbound}, \eqref{eq:parabolicupperbound}, and using \eqref{eq:final eq} for $r=l_1$, we obtain the desired estimates. 
			
		\end{proof}
	\end{lemma}

	\begin{proof}[\bf Proof of Theorem \ref{thm:extremal_length}]
		Let $h$ denote the associated Koenigs function of $(\phi_t)$. 		
		
		\textbf{(i) Hyperbolic Case.} Suppose that the petal $\Delta$ has $\alpha$-point $\sigma \in \partial \D$, different from the Denjoy--Wolff point of $(\phi_t)$, and is thus hyperbolic. We handle the two cases of non-elliptic and elliptic semigroups individually. 
		
		Assume, first, that $(\phi_t)$ is a non-elliptic semigroup. In this case, $h(\Delta)$ is a maximal horizontal strip $S$ inside the Koenigs domain of width equal to $-\frac{\pi}{\mu}.$
		The image of $K$ under $h$ is a compact set in $\overline{S}$. Set $l_1: = \min \{ \RE w: w \in h(K) \} $, $l_2 : = \max \{\RE w : w \in h(K) \}$, and $\theta_x:=\overline{S}\cap\{w \in \mathbb{C} : \textup{Re} w=x\}$, for $x\in \mathbb{R}$. Each $\theta_x$ is a vertical line segment in $\overline{S}$ with endpoints on $\partial S$. Hence, each $\theta_x$ is a crosscut of $S$. 
		Suppose that $l_1<l_2$. Fix $t_0>0$ such that $l_2-t_0 <l_1$. Note that the segments $\theta_{l_1}$ and $\theta_{l_2 -t}$ separate $h(K)$ and $h(\phi_{-t}(K))= h(K)-t$ inside $S$, for all $t\geq t_0$. 
		Let us denote by $\Gamma_t^{\star}$ the family of curves that separate $h(K)$ from $h(K)-t$ in $S$. Set $\Gamma_t^{1\star}$ to be the family of curves in $S$ separating $\theta_{l_1}$ and $\theta_{l_2 -t}$. 
		Obviously, $\Gamma_t^{1\star} \subset \Gamma_t^{\star}$ and thus, by the extension rule (Lemma \ref{lem:extension}) of the extremal length and the formula of extremal distance for rectangles \eqref{eq:rectangle}, we obtain 
		\begin{equation}\label{eq:lowerbound}
			\lambda (\Gamma_t^{\star}) \leq \lambda (\Gamma_t^{1\star})= \lambda_S^{\star}(\theta_{l_1}, \theta_{l_2-t})=-\frac{\pi}{\mu (l_1-l_2+t)}. 
		\end{equation}
		Hence $\lambda_{S}(h(K), h(K)-t) \geq  -\frac{\mu t }{\pi} -\mu \frac{l_1-l_2}{\pi}$. 
		We now focus on obtaining an upper bound for the extremal distance $\lambda_S (h(K), h(K)-t)$. Set $\Gamma_t^{2}$ the family of curves in $S$ that join $\theta_{l_1-t}$ and $\theta_{l_2}$ and $\Gamma_t^{2\star}$ the family of curves in $S$ that separate $\theta_{l_1-t}$ from $\theta_{l_2}$. Note that the latter consists of all crosscuts of $S$ in $S\cap\{w: \text{Re} w\in [l_1-t, l_2]\}$. 
		
		There exists a subfamily $\Gamma_{C}^{\star}$ of $\Gamma_t^{\star}$ that contains crosscuts of $S$ that simultaneously separate $\theta_{l_1-t}$ from $\theta_{l_2}$. Thus, by the extension rule of the extremal length, 
		\begin{equation}\label{eq:crosscuts}
			\lambda(\Gamma_{C}^{\star}) \geq \lambda (\Gamma_t^{2\star})= \frac{1}{\lambda_S (\theta_{l_1-t}, \theta_{l_2})}.
		\end{equation}
		
		The rest of the curves of $\Gamma_t^{\star}$ fall in the following categories:
		\begin{enumerate}
			\item closed curves $\gamma $ with trace in $S$, where $h(K) \subset \overline{\textup{int}(\gamma)}$;
			\item closed curves $\gamma $ with trace in $S$, where $h(K)-t \subset \overline{\textup{int}(\gamma)}$; 
			\item crosscuts with endpoint(s) on $\partial S$ that intersect $\theta_{l_1}$ (resp. $\theta_{l_2-t}$) and $\theta_{l_2}$ (resp. $\theta_{l_1-t}$).
		\end{enumerate}
		
		Denote by $F_1,F_2$, and $F_3$, respectively, the three families above. 
		In the families $F_1$ and $F_2$, the closed curves intersect $\theta_{l_1}$ and $\theta_{l_2}$ or $\theta_{l_1-t}$ and $\theta_{l_2-t}$, respectively. Note that these cases do not appear when $K\cap \partial \Delta \neq \emptyset$. 
		As a result, for all the respective families of curves in the aforementioned three cases, it holds that every curve $\gamma^{\prime}$ contains some curve that joins either $\theta_{l_1}$ with $\theta_{l_2}$ or $\theta_{l_1-t}$ with $\theta_{l_2-t}$. Thus, for the extremal length of each family we have 
		\begin{equation}\label{eq:extremal length maximum}
			\lambda(F_j)\ge\max\{ \lambda_S(\theta_{l_1}, \theta_{l_2}), \lambda_S(\theta_{l_1-t}, \theta_{l_2-t}) \} =\frac{1}{\lambda^{\star}_S(\theta_{l_1},\theta_{l_2})}= - \frac{\mu(l_2-l_1)}{\pi}, \quad\text{for }j=1,2,3.  
		\end{equation}
		Following Lemma \ref{non-disjointparallel} and applying \eqref{eq:conjugate}, we get
		\begin{eqnarray}\label{eq:upperbound}
			\lambda(\Gamma_t)&=\frac{1}{\lambda(\Gamma_t^{*})}\leq& \frac{1}{\lambda(\Gamma_C^{\star})}+ \sum_{i=1}^3 \frac{1}{\lambda(F_i)}\nonumber \\
			&\underset{\eqref{eq:crosscuts},\eqref{eq:extremal length maximum}}{\leq}& 
			\lambda_S (\theta_{l_1-t}, \theta_{l_2}) - \frac{3\pi}{ \mu (l_2-l_1)}  \nonumber \\ 
			&=& -\frac{\mu t}{\pi} - \left( \mu \frac{l_2-l_1}{\pi} + \frac{3\pi}{\mu (l_2 -l_1)} \right).
		\end{eqnarray}
		Combining \eqref{eq:lowerbound} and \eqref{eq:upperbound}, we get $t^{-1}\lambda(\Gamma_t) \to -\mu/\pi$, as $t\to +\infty$. 
		
		Let us now consider the case where $l_1=l_2$. This implies that $K$ is a vertical line segment and $h(K) \subseteq \theta_{l_1}$.
		The lower bound for $\lambda_{\Delta}(K, \phi_{-t}(K))$ follows again from \eqref{eq:lowerbound}. In this case, $\Gamma_C^{\star}$ contains all the curves separating $h(K)$ from $h(K)-t$ that simultaneously separate $\theta_{l_1}$ from $\theta_{l_1-t}$. The family of curves $F_3$ consists of all crosscuts with endpoint(s) in $\partial S$ that intersect $\theta_{l_1} \setminus h(K)$ (or resp. $\theta_{l_1-t}\setminus (h(K)-t)$, for $t$ large enough). Set $m_1:=\min \{ \textup{Im} w: w \in h(K)\}$ and $m_2 := \max \{ \textup{Im} w: w \in h(K)\}$. Every curve of the families $F_j$, $j=1,2,3$, either joins the horizontal lines $\{w \in \mathbb{C} : \textup{Im}w= m_i\} $, $i=1,2$, inside the rectangle $[l_1-1, l_1] \times [m_1, m_2]$ (resp. $[l_1-t, l_1+1-t] \times [m_1, m_2]$) or joins $\theta_{l_1}$ with $\theta_{l_1-1}$ (resp. $\theta_{l_1-t}$ with $\theta_{l_1+1-t}$). Thus, we obtain 
		$$\frac{1}{\lambda(\Gamma_C^{\star})}\le-\frac{\mu t}{\pi}\quad \text{and}\quad\frac{1}{\lambda(F_j)}\le \frac{1}{\textup{diam}(h(K))} -\frac{\pi}{ \mu} \, \quad  j=1,2,3. $$ 
		Hence, using Lemma \ref{non-disjointparallel}, there exists a constant $C=C(K)$ such that in \eqref{eq:upperbound} it holds $\lambda_{\Delta}(K, \phi_{-t}(K)) \leq -\frac{\mu}{\pi} t+ C$, which leads to the desired result.

		We now consider the case where $(\phi_t)$ is an elliptic semigroup. Without loss of generality, we may suppose that $\tau=0$ is the Denjoy--Wolff point of $(\phi_t)$ with spectral value $\lambda\in \mathbb{C}$, $\mathrm{Re}\lambda>0$. The image of $\Delta$ under $h$ is a maximal $\lambda$-spirallike sector of center $e^{i \theta_0}$ for some $\theta_0 \in (-\pi, \pi]$ and amplitude $2 a:=-\frac{|\lambda|^2 \pi}{\mu \text{Re} \lambda}$.
		With the aid of Lemma \ref{lem:mapofpetal}, we use the conformal mapping $g(z) = \frac{e^{-i \Arg \lambda}}{\cos(\Arg \lambda)} \Log z + i (a-\theta_0)$, which maps $h(\Delta)$ onto the horizontal strip $S_a:=\{z :  0< \text{Im} z <2a\}$.
		Let $w\in h(\Delta)$. By simple calculations, $g(e^{\lambda t} w)= t \frac{|\lambda|}{\cos(\Arg \lambda)} + g(w)$. Set $c(\lambda)\vcentcolon=\frac{|\lambda|}{\cos(\Arg \lambda)} > 0$ and $\widehat{K}\vcentcolon =g(h(K))$. Under this notation, $g(e^{\lambda t} h(K)) =  \widehat{K} + t c( \lambda)$. Let us write $l_1\vcentcolon = \min \{ \RE w: w \in \widehat{K}\} $ and $l_2 \vcentcolon = \max \{\RE w : w \in \widehat{K} \}$. Suppose $l_1< l_2$. We follow the same procedure as in the non-elliptic case and we obtain the estimates
		$$\frac{l_1 -l_2}{2a}+ \frac{c(\lambda) t}{2a} \leq \lambda_{S_a}(\widehat{K}, \widehat{K}+c(\lambda) t ) \leq   \frac{c(\lambda) t}{2a} + \left( \frac{l_2-l_1}{2a} + \frac{6a}{l_2 -l_1} \right). $$
		Therefore
		$$\lim_{t\to +\infty} \frac{\lambda_{S_a}(\widehat{K}, \widehat{K}+c(\lambda) t )}{t} = \frac{c(\lambda)}{2a} =  -\frac{|\lambda|\mu \text{Re}(\lambda)}{|\lambda|^2 \pi \cos(-\Arg \lambda)} =  -\frac{\mu}{\pi} $$
		and due to the conformal invariance of the extremal length, we obtain the desired result.  If $l_1=l_2$, then $h(K)$ is a spirallike segment in $h(\Delta)$, and following the same arguments as in the hyperbolic case where $h(K)$ is a vertical slit, we again obtain the desired result. 
		
		\textbf{(ii) Parabolic Case.} Set $H:=h(\Delta)$ and suppose it is of the form $\{z \in \mathbb{C} : \text{Im} z> a\}$, for some $a\in \mathbb{R}$. The proof for the case where $h(\Delta)$ is a lower half-plane follows similarly. Due to the conformal invariance of the extremal length, we have 
		$$\lambda_{\Delta}(K, \phi_{-t}(K)) = \lambda_H(h(K), h(K)-t ) =\lambda_{\mathbb{H}}(h(K)-i a, h(K)-i a -t),$$
		where $\mathbb{H}$ denotes the usual upper half-plane. By setting $E:=h(K)-i a$, we will apply Lemma \ref{lem:extremallengthcircles} to obtain bounds for $\frac{\lambda_{\Delta}(K, \phi_{-t}(K)) }{\log t}$. 
		Suppose first that $i a \notin h(K)$ and that $h(K)$ is not a circular arc of center $i a$. 
		Set $l_1:=\min \{ |w| :w\in E\} = \min\{|w -ia| : w\in h(K)\} $ and $l_2:=\max \{ |w| :w\in E\} $. Then $l_1<l_2$ and according to Lemma \ref{lem:extremallengthcircles}, for all $t>2l_2$, 
		\begin{equation}\label{eq:inequality_parabolic}
			\frac{2}{\pi \log t} \log \frac{1+\sqrt{1-4\frac{l_2^2}{t^2}}}{2l_2} \leq \frac{\lambda_{\Delta}(K, \phi_{-t}(K))}{\log t}- \frac{2}{\pi}  \leq \frac{2}{\pi\log t } \log \frac{1+\sqrt{1-4\frac{l_1^2}{t^2}}}{2l_1} +  \frac{2 \log \frac{l_2}{l_1} }{\pi\log t} 
		\end{equation}  
		leading to 
		$$\lim_{t\to +\infty} \frac{\lambda_{\Delta}(K, \phi_{-t}(K)}{\log t} =\frac{2}{\pi} .$$
		In the case $i a \in h(K)$ or $h(K)$ is a circular arc centered at $i a$ (we mention this case separately because it either implies $l_1=0$ or $l_1=l_2$, respectively), we set $p \in \mathbb{R}^{+}$ such that $p=2 \max\{ |w-ia| : w\in h(K)\} $ and $l_1:=\min\{|w-p -ia| : w\in h(K)\}$ and $l_2:=\min \{ |w-p -i \alpha| :w\in h(K)\} $. Then $E \subset \overline{A(p, l_1,l_2)}$ and since the extremal length remains invariant under translations, we may consider $E^{\prime} = E-p$ and $\lambda_{\mathbb{H}}(E^{\prime}, E^{\prime}-t)=\lambda_{\Delta}(K, \phi_{-t}(K))$. We apply Lemma \ref{lem:extremallengthcircles}, and we derive again \eqref{eq:inequality_parabolic}. 
	\end{proof}

	\begin{remark}
		As also referred to in the Introduction, the extremal distance can be viewed as the Dirichlet integral of the solution to the Dirichlet--Neumann problem (\cite[Theorem 4-5]{ahlforsconformal}). Indeed, via Theorem \ref{thm:extremal_length}, we obtain information on the asymptotic behavior
		of the capacity of the generalised condenser $(\Delta,\{K,\phi_{-t}(K)\},\{0,1\})$ , as $t\to +\infty$. For further information on condensers, the interested reader may refer to \cite{dubininbook}.     
	\end{remark}
	
	Observe that Theorem \ref{thm:extremal_length} allows the compact set $K$ to fully lie on the boundary of a petal $\Delta$. In such a case, we say that the backward orbit of $K$ is non-regular (this subject will be further explored in the next section). However, we see that regularity does not impact the divergence of the extremal length when the backward orbit shrinks to the Denjoy--Wolff point or a repelling fixed point of the semigroup. This behavior is not shared by the hyperbolic length as we shall see later on. To end the section, we proceed to an extremal example, where the compact set whose rate we examine, indeed lies on the boundary of a petal.
	\begin{example}
		Let $(\phi_t)$ be a parabolic semigroup with associated Koenigs function $h$. Let $\Delta$ be a parabolic petal. The image $h(\Delta)$ is a horizontal half-plane. Without loss of generality, suppose $h(\Delta)=\{ w\in \mathbb{C} : \textup{Im} w> 0 \}$. Suppose $h(K)=[1,2]$ and $h(\phi_{-t}(K))=[1-t,2-t]  $. Using a conformal automorphism $T_t$ of $\mathbb{H}$ with $T_t(w)= \frac{w-(2-t)}{w-2}$, for $t>2$, we map the intervals $[1-t,2-t]$ and $[1,2]$ onto the intervals $[0,\frac{1}{t+1}]$ and $(-\infty,1-t]$, respectively. Set $\lambda_t:=\lambda_{\mathbb{H}}(T_t([1,2]),T_t([1-t,2-t]))$. According to \cite[Ex. 4(c), p. 151]{margarnett}, it follows that 
		$$\frac{8}{\pi} e^{-\pi \lambda_t/2} \sim \frac{2}{\pi} \arctan\sqrt{\frac{\frac{1}{t+1}}{t-1}}=\frac{2}{\pi} \arctan\sqrt{ \frac{ 1}{t^2-1}}$$
		and thus, $$-\frac{\pi \lambda_t}{2} \sim -2\log 2+ \log \arctan \sqrt{ \frac{1}{t^2-1}} ,$$
		where the notation $\sim$ signifies the similar asymptotic behavior, as $t\to+\infty$. Due to the conformal invariance of extremal length, we are led through elementary calculations to 
		$$\lim_{t\to +\infty } \frac{\lambda_{\Delta}(K, \phi_{-t}(K))}{\log t} = \frac{2}{\pi} \lim_{t\to +\infty}\left( -\frac{\log \arctan \sqrt{ \frac{1}{t^2-1}} }{\log t}\right) = \frac{2}{\pi} .$$
	\end{example}

	\section{Non-Regular Orbits}\label{sec:non-regular}

	As stated in Subsection \ref{sub:semigroups}, a non-regular backward orbit $\tilde{\gamma}_z:[0,+\infty)\to\mathbb{D}$ with starting point $z\in \mathcal{W}$ converges to a boundary fixed point and its image under the associated Koenigs function is either a horizontal half-line $\{h(\tilde{\gamma}_z(0))-t:t\geq 0\}$ or a half-spiral $\{e^{\lambda t} h(\tilde{\gamma}_z(0)) : t \geq 0 \}$. Non-regular orbits are more volatile and present more intricacies in their behavior in comparison to their regular counterparts. In this final section of the article, we are going to investigate their divergence rates.

	\subsection{Hyperbolic length}
	We commence with the hyperbolic length and we prove a first result bounding the corresponding asymptotic behavior of non-regular backward orbits.

	\begin{proposition}\label{prop:non-reg}
		Let $(\phi_t)$ be a non-elliptic semigroup in $\D$ with Koenigs function $h$ and Koenigs domain $\Omega$. Suppose that for some $z\in\D$ the backward orbit $\tilde{\gamma}_z:[0, +\infty) \to \D$ is well-defined and non-regular (that is $z\in\D\cap\partial\mathcal{W}$). Then $$\lim_{t\to +\infty} \frac{\ell_{\D}(\tilde{\gamma}_z;[0,t])}{t} =+ \infty\, , $$ whereas $$\limsup_{t\to +\infty} \frac{\ell_{\D}(\tilde{\gamma}_z;[0,t])}{\frac{t}{\delta_{\Omega}(h(z)-t)}} \leq 1\, .$$
		\begin{proof}
			We use the connection between the hyperbolic density and the Euclidean distance to the boundary. Fix $t>0$. By the conformal invariance of the hyperbolic length, we have $$\ell_{\D}(\tilde{\gamma}_z; [0,t])= \ell_{\Omega}(h \circ \tilde{\gamma}_z ;[0,t]) = \int_0^t \lambda_{\Omega}(h(z)-s) ds \geq \frac{1}{4} \int_0^t \frac{ds}{\delta_{\Omega}(h(z)-s)} \, . $$
			The function $s\mapsto \delta_{\Omega}(h(z)-s)$ is (not necessarily strictly) decreasing to $0$ in $[0, +\infty)$, due to the geometry of $\Omega$. Thus, given an arbitrarily small $\epsilon>0$, we may find $t_0 \geq 0$ such that $\delta_{\Omega} (h(z)-t) < \epsilon$, for all $t \geq t_0$. Suppose $t \geq t_0$ and set $c:=\frac{1}{4} \int_0^{t_0} \frac{ds}{\delta_{\Omega}(h(z)-s)}$. Following the above calculations, we have 
			$$\ell_{\D}(\tilde{\gamma}_z; [0,t]) \geq  c +\frac{1}{4} \int_{t_0}^t \frac{ds}{\epsilon} = c+ \frac{t-t_0}{4\epsilon}$$ 
			and therefore, $$\liminf_{t\to +\infty} \frac{\ell_{\D}(\tilde{\gamma}_z; [0,t])}{t} \geq \frac{1}{4\epsilon}.$$
			The above is true for any choice of $\epsilon$ and letting $\epsilon \to 0$, we obtain that the limit $\frac{\ell_{\D}(\tilde{\gamma}_z; [0,t])}{t}  \to +\infty$, as $t \to +\infty$. 
			With a similar series of ideas, we obtain that $\ell_{\D}(\tilde{\gamma}_z; [0,t]) \leq \int_0^t \frac{ds}{\delta_{\Omega}(h(z)-s)} $. However, since $\delta_{\Omega}(h(z)-s) \geq \delta_{\Omega }(h(z) -t)$ for all $s\in [0,t]$, we obtain 
			$$\ell_{\D}(\tilde{\gamma}_z; [0,t]) \leq \frac{1}{\delta_{\Omega}(h(z)-t)} \int_0^t ds=\frac{t}{\delta_{\Omega}(h(z)-t)}$$
			and hence, we conclude the desired inequality. 
			
		\end{proof}
	\end{proposition}
	
	The preceding proposition indicates that for non-regular backward orbits the hyperbolic length diverges to $+\infty$ much faster than $t$ does. We should also note that if the hyperbolic length diverges to $+\infty$ slower than $\frac{t}{\delta_{\Omega}(h(z)-t)}$, where $z$ is the starting point of the non-regular backward orbit, then we cannot derive further information. In order to verify the sharpness of the limiting behavior above, we state the following example. 
	
	\begin{example}
		\normalfont
		Suppose $(\phi_t)$ is a non-elliptic semigroup in $\D$ with Koenigs function $h$ and Koenigs domain $\Omega$. Suppose that for some $z\in \D$, we have $$\delta_{\Omega}(h(z)-t)= \begin{cases} \frac{1}{\ln 2}, & t \in [0,2] \\
			\frac{1}{\ln t} , & t \in [2, +\infty)            .
		\end{cases}$$
		Note that $t\mapsto \delta_{\Omega}(h(z)-t)$ is continuous, decreasing and satisfies $\delta_{\Omega}(h(z)-t)\to 0$, as $t \to +\infty$, exactly as non-regular backward orbits behave; see \eqref{distancedegen}. For $t>2$, it is trivial to see that $$\int_0^t \frac{ds}{\delta_{\Omega}(h(z)-s)} = t \ln t - t+2. $$
		Using \eqref{eq:relationquasi}, we get $$\frac{ t \ln t - t+2}{4} \leq \ell_{\D}(\tilde{\gamma}_z ; [0,t]) \leq  t \ln t - t+2, $$ for all $t\ge 2$, leading to $$\frac{1}{4} \leq \liminf_{t \to +\infty} \frac{\ell_{\D}(\tilde{\gamma}_z ; [0,t])}{\frac{t}{\delta_{\Omega}(h(z)-t)}} \leq \limsup_{t \to +\infty} \frac{\ell_{\D}(\tilde{\gamma}_z ; [0,t])}{\frac{t}{\delta_{\Omega}(h(z)-t)}} \le 1,$$
		which renders the upper estimate in Proposition \ref{prop:non-reg} sharp. 
	\end{example}
	
	When it comes to the sharpness of the lower bound of the rate by which the hyperbolic length along non-regular backward orbits diverges to $+\infty$, we prove the following proposition. 
	
	\begin{proposition}
		Let $f:[0,+\infty) \to [0,+\infty) $ be a convex, differentiable function such that 
		$$\lim_{t\to+\infty} \frac{f(t)}{t} = +\infty.$$ 
		Then, there exist a non-elliptic semigroup $(\phi_t)$ in $\D$ and a point $z \in \D$ such that $\frac{\ell_{\D}(\tilde{\gamma}_z ; [0,t])}{f(t)} = O(1)$, as $t\to +\infty$. 
		\begin{proof}
			Due to the hypothesis $f^{\prime}$ is well-defined. Since $f$ is convex, $f^{\prime}$ is increasing and the graph of $f$ lies above all of its tangents; namely $f(0) \geq f(t)-f^{\prime}(t) t$, for all $t\geq 0$. Therefore $f^{\prime}(t) \ge \frac{f(t)-f(0)}{t}$, for all $t>0$, leading to $ f^{\prime}(t) \to +\infty$, as $t\to +\infty$. All in all, the function $\frac{1}{f^{\prime}}$ is well-defined in an interval of the form $[a, +\infty)$, where $a\geq 0$, and decreases to $0$. Without loss of generality, assume that $a=0$. We may construct a simply connected domain $\Omega$, convex in the positive direction with $w \in \Omega$ such that $\delta_{\Omega}(w-t) =\frac{1}{f^{\prime}(t)}$, for all $t \geq 0$. This time we do not need a condition akin to that in \eqref{eq:condition} because in the backward setting, the distance $\delta_{\Omega}(w-t)$ is allowed to decay to $0$ uncontrollably fast. Let $h:\D\to \Omega$ be a Riemann map of $\Omega$ and set $z:=h^{-1}(w)$. The relation $\phi_t(\cdot) : = h^{-1}(h(\cdot)+t)$ defines a non-elliptic semigroup in $\D$ and by construction $\tilde{\gamma}_z$ is a non-regular backward orbit of $(\phi_t)$. Combining the above arguments with \eqref{distancedegen}, we obtain
			$$\frac{1}{4} \int_0^t f^{\prime} (s)ds  \leq \ell_{\D}(\tilde{\gamma}_z ; [0,t]) \leq \int_0^t f^{\prime}(s) ds $$ and thus,
			since $\int_0^t f^{\prime} (s)ds = f(t)-f(0)$, we conclude 
			$$\frac{1}{4} \leq \liminf_{t\to +\infty} \frac{\ell_{\D}(\tilde{\gamma}_z ; [0,t]) }{f(t)} \leq \limsup_{t\to +\infty}  \frac{\ell_{\D}(\tilde{\gamma}_z ; [0,t]) }{f(t)} \leq 1,$$
			which yields the desired outcome. 
		\end{proof}
	\end{proposition}
	As a result, we understand that the rhythm at which the hyperbolic length diverges is always faster than linear. However, there is no other rhythm which is simultaneously faster than linear, but slower than the divergence of the hyperbolic length of every non-elliptic semigroup. This provides a sense of sharpness for the lower bound in Proposition \ref{prop:non-reg}.
	
	We move on to elliptic semigroups. The proof of the following result resembles the ones executed before. For this reason, we opt to only sketch it. 
	\begin{proposition}\label{prop:elliptic non-reg}
		Let $(\phi_t)$ be an elliptic semigroup in $\D$ with spectral value $\lambda\in\mathbb{C}$, $\mathrm{Re}\lambda>0$. Suppose that for some $z\in\D$ the backward orbit $\tilde{\gamma}_z:[0,+\infty)\to\D$ is well-defined and non-regular (that is $z\in\D\cap\partial\mathcal{W}$). Then
		\begin{equation*}
			\lim_{t\to+\infty}\frac{\ell_{\D}(\tilde{\gamma}_z;[0,t])}{e^{\mathrm{Re}\lambda t}}=+\infty.
		\end{equation*}
	\end{proposition}
	\begin{proof}
		Fix $t>0$. Recall that $h(\tilde{\gamma}_z(t))=e^{\lambda t}h(z)$ since we deal with an elliptic semigroup. Following the steps introduced in previous proofs, we have
		\begin{equation}\label{eq:ell, non-reg 1}
			\ell_{\D}(\tilde{\gamma}_z;[0,t])\ge\frac{\lvert\lambda h(z)\rvert}{4}\int_{0}^{t}\frac{e^{\mathrm{Re}\lambda s}ds}{\delta_{\Omega}(e^{\lambda s}h(z))}.
		\end{equation}
		Once more, $\delta_{\Omega}(e^{\lambda s}h(z))\to0$, as $t\to+\infty$, since $\tilde{\gamma}_z$ is non-regular. As a result, given $\epsilon>0$, we may find $t_0\ge0$, exactly as in preceding proofs, so that \eqref{eq:ell, non-reg 1} leads to
		\begin{equation*}
			\ell_{\D}(\tilde{\gamma}_z;[0,t])\ge c+\frac{\lvert\lambda h(z)\rvert}{4\epsilon}\int_{t_0}^{t}e^{\mathrm{Re}\lambda s}ds=c+\frac{\lvert\lambda h(z)\rvert}{4\epsilon}\left(e^{\mathrm{Re}\lambda t}-e^{\mathrm{Re}\lambda t_0}\right).
		\end{equation*}
		Dividing by $e^{\mathrm{Re}\lambda t}$, taking limits as $t\to+\infty$, and finally letting $\epsilon \to 0$, we reach the desired outcome.
	\end{proof}
	It is straightforward that a combination of Theorem \ref{thm:hyplength_backward} and Propositions \ref{prop:non-reg}, \ref{prop:elliptic non-reg} directly yields Corollary \ref{cor:non-reg} which serves as a counterpart of Corollary \ref{cor:hyperbolic step} for backward orbits.

	\subsection{Extremal length}
	We proceed to the analogous investigation concerning extremal length. Nevertheless, before we move on to backward orbits of compact sets, we will need the following information from \cite{BCZZ_25} on non-regular backward orbits converging to super-repelling fixed points of semigroups. 
	
	\begin{lemma}\cite[Theorem 2.2 (d)]{BCZZ_25}\label{lem:Rates}
		Let $(\phi_t)$ be a non-elliptic semigroup in $\mathbb{D}$ having a non-regular backward orbit $\widetilde{\gamma}_z:[0,+\infty) \to \D$, $z\in \D$, which lands at a super-repelling fixed point $\xi \in \partial \D$.  Then for every $\epsilon>0 $, there exists $C\vcentcolon=C(z, \epsilon)>0$ such that 
		$$|\widetilde{\gamma}_z(t)-\xi|\leq C e^{-\epsilon t}.$$
	\end{lemma}

	Let us now consider a compact, non-polar set $K$ that lies in the backward invariant set $\mathcal{W}$ of a semigroup $(\phi_t)$ such that for all $z\in K$, the backward orbits $\tilde{\gamma}_z$ converge to the same super-repelling fixed point of $(\phi_t)$. At once, this implies $K\subset\partial\mathcal{W}$. Without loss of generality, we assume that $K$ is a continuum. Hence $h(K)$ is either a line segment (non-elliptic case) or a spiral arc (elliptic case). 
	
	As stated in \eqref{distancedegen} and denoting by $h$ the corresponding Koenigs function, an important property is that all points $z$ of $K$ satisfy $\delta_{\Omega}(h(\phi_t(z))=\textup{dist}(h(\phi_t(z)), \partial \Omega)\to0$, as $t\to +\infty$. Through the work executed below, we are going to discover that the distance of the orbit's image under the Koenigs function from the boundary of the Koenigs domain plays a profound role in the investigation of the respective extremal length. Towards this goal, for $w\in\Omega$, we define
	$$\delta_{\Omega}^+(w) \vcentcolon= \inf\left\{\textup{Im} p -\textup{Im} w : p \in \partial \Omega\, , \, \, \textup{Re} p =\mathrm{Re}w \, \,  \text{and} \, \,  \textup{Im} p > \textup{Im} w \right\}$$
	and 
	$$\delta_{\Omega}^-(w)\vcentcolon = \inf\left\{\textup{Im} w - \textup{Im} p :  p \in \partial \Omega\, , \, \, \textup{Re} p =\mathrm{Re}w \, \,  \text{and} \, \,  \textup{Im} p < \textup{Im} w \right\}$$
	to be the \textit{vertical distances} of $w$ from the boundary of $\Omega$. 
	Let $F$ be a compact set in $\Omega$. Set $\delta_{\Omega}^+(F) = \min_{w \in F} \delta_{\Omega}^+ (w)$ and $\delta_{\Omega}^-(F) = \min_{w \in F} \delta_{\Omega}^- (w)$. 
	
	We obtain the following result on the asymptotic behavior, as $t\to+\infty$, of the extremal distance $\lambda_{\D}(K, \phi_{-t}(K))$, in terms of the Euclidean geometry of the associated Koenigs domain of the semigroup.

	\begin{proposition}\label{prop:non-regular}
		Let $(\phi_t)$ be a non-elliptic semigroup in $\D$ with associated Koenigs domain $\Omega$ having a non-regular backward orbit $\widetilde{\gamma}_z:[0,+\infty) \to \D$, $z\in \D$, which lands at a super-repelling fixed point $\xi \in \partial \D$. Let $K$ be a compact non-polar subset of $\tilde{\gamma}_z([0,+\infty))$ and set $\widetilde{K}:=h(K)$. Then 	
		$$\limsup_{t\to +\infty} \left[ \lambda_{\Omega}(\widetilde{K}, \widetilde{K}-t) \cdot \frac{\delta_{\Omega}^+(\widetilde{K}-t)+ \delta_{\Omega}^-(\widetilde{K}-t)}{\ell_e(\widetilde{K})+t} \right] \leq 1 \, , $$	
		where $\ell_e$ denotes the Euclidean length. 
		\begin{proof}
			Suppose $a:=\textup{Im}\,  h(z)$ something which means that $\widetilde{K} \subset i a+\mathbb{R}$. Set $l_1:=\min \{ \textup{Re} w: w \in \widetilde{K}\}$ and $l_2:=\max \{ \textup{Re} w: w \in \widetilde{K}\}$.
			Note that due to the property of $\Omega$ being convex in the positive direction, we infer
			$$h_t:=\delta_{\Omega}^+(\widetilde{K}-t)+ \delta_{\Omega}^-(\widetilde{K}-t)=\delta_{\Omega}^+ (l_1-t+ i a)+\delta_{\Omega}^- (l_1-t+ i a).$$ 
			Then the horizontal rectangle 
			$$S_t:= \left\{ x+iy \in \mathbb{C} : x \in (l_1-t, l_2) \quad \& \quad y\in (a-\delta_{\Omega}^- (l_1-t+ia),a+ \delta_{\Omega}^+ (l_1-t+ia))\right\}$$
			of width $l_2-l_1+t$ and of height $h_t$ is contained inside $\Omega$. Therefore, by the extension rule of extremal length (Lemma \ref{lem:extension}), we obtain
			\begin{equation}\label{eq:ineq_degeneratestrip}
				\lambda_{\Omega}(\widetilde{K}, \widetilde{K}-t) \leq \lambda_{S_t}(\widetilde{K}, \widetilde{K}-t). 
			\end{equation} 
			Taking the properties of the extremal length in rectangles into consideration, we obtain 
			$$ \lambda_{S_t}(\widetilde{K}, \widetilde{K}-t) \geq  \frac{l_2-l_1}{h_t} +\frac{t}{h_t} $$ 
			and $$ \lambda_{S_t}(\widetilde{K}, \widetilde{K}-t) \leq   \frac{l_2-l_1}{h_t} +\frac{t}{h_t} + \frac{3h_t}{l_2-l_1},$$
			where the second inequality follows by applying the methodology of the proof in the hyperbolic case of Theorem \ref{thm:extremal_length}. As a result, $$\lambda_{S_t}(\widetilde{K}, \widetilde{K}-t)- \left( \frac{l_2-l_1}{h_t} +\frac{t}{h_t}\right) \in \left[0, \frac{3 h_t}{l_2-l_1} \right]$$
			and due to \eqref{distancedegen} and \cite[Corollary 13.6.7]{Booksem}, which states that the unrestricted limit of $\textup{Im} h(w)$ goes to $\textup{Im} h(\xi)$ as $w\to \xi$ , we have 
			$$\lim_{t\to +\infty}\frac{\lambda_{S_t}(\widetilde{K}, \widetilde{K}-t)}{\frac{t+ l_2-l_1}{h_t} } = 1$$ which combined with \eqref{eq:ineq_degeneratestrip} leads to the desired result. 		 
		\end{proof}
	\end{proposition}
	
	Let us observe the following example showing how fast the extremal distance $\lambda_{\D}(K, \phi_t(K))$ can diverge to $+\infty$.
	
	\begin{example}
		Let $(\phi_t)$ be a parabolic semigroup of zero hyperbolic step with associated Koenigs domain $$\Omega=\{x+iy \in \mathbb{C} :  x\in \mathbb{R}, \, -e^x < y< e^x\}.$$
		Set $\widetilde{K}=[1,2]\subset \mathbb{R}$. 
		Then $h_t=\delta_{\Omega}^+(\widetilde{K}-t)+ \delta_{\Omega}^-(\widetilde{K}-t)= 2 e^{1-t}$. Applying Proposition \ref{prop:non-regular}, we get $$\limsup_{t\to +\infty}\lambda_{\Omega} (\widetilde{K}, \widetilde{K}-t)\frac{2e^{1-t}}{t+1} \leq 1.$$As a result, $\lambda_{\Omega} (\widetilde{K}, \widetilde{K}-t)$ diverges to $+\infty$ at most by rate $te^t$.  
	\end{example}

	Moreover, we obtain the following asymptotic behavior of extremal distance. 
	
	\begin{proposition}
		Let $(\phi_t)$ be a semigroup in $\D$ having a non-regular backward orbit $\widetilde{\gamma}_z:[0,+\infty) \to \D$, $z\in \D$, which lands at a super-repelling fixed point $\xi \in \partial \D$. Let $K$ be a non-polar compact subset of $\tilde{\gamma}_z([0,+\infty))$. Then 
		$$\lim_{t\to +\infty} \frac{\lambda_{\D}(K, \phi_{-t}(K))}{t}=+\infty.$$
		\begin{proof}
			Without loss of generality, assume that $K$ is connected and thus, there exist $a, b\in \mathbb{R}$ with $a<b$ such that $K=\{\phi_{-x}(z) : x\in [a,b] \}$. Set $l:=\min \{|\xi-\phi_{-x}(z)| : x\in [a,b]\}=\text{dist}(\xi,K)$. Pick $t>0$ sufficiently large so that $l> \max\{ |\xi-\phi_{-t}(w)| : w \in K\}$ and find $t_{a,b}$ such that 
			$$m_{t}:=|\phi_{-t_{a,b}}(z)-\xi| = \max_{t^{\prime} \in [t+a, t+b]} |\phi_{-t^{\prime}}(z)-\xi|.$$  
			Let $H$ be the half-plane that is tangent to the unit disk at $\xi$ and contains the unit disk. Fix an arbitrary $\epsilon>0$. By Lemma \ref{lem:Rates}, there exists a positive constant $C\vcentcolon=C(z,\epsilon)$ so that
			\begin{equation}\label{eq:final 1}
				|\xi-\phi_{-t_{a,b}}(z)|\le Ce^{-\epsilon t_{a,b}}.
			\end{equation}
			Then, for the extremal length we have
			\begin{eqnarray*}
				\lambda_{\D}(K, \phi_{-t}(K)) &\geq &\lambda_{\D}(\partial D(\xi, l)\cap \D, \partial D(\xi, m_t)\cap \D)\\
				&\geq & \lambda_{H}(\partial D(\xi, l)\cap H, \partial D(\xi, m_t)\cap H)\\
				&\underset{\text{Lemma }\ref{lem:halfannulus}}{=}&\frac{1}{\pi} \log \frac{\text{dist}(\xi, K)}{ |\xi- \phi_{-t_{a,b}}(z_0)|}\\
				&\underset{\eqref{eq:final 1}}{\geq}& \frac{1}{\pi} \log \text{dist}(\xi, K) -\frac{1}{\pi} \log C -\frac{1}{\pi} \log e^{-\epsilon t_{a,b}} \\
				&\ge& \frac{1}{\pi} \log \text{dist}(\xi, K) -\frac{1}{\pi} \log C +\frac{\epsilon}{\pi} (t+a).
			\end{eqnarray*}
			As a result,dividing $t>0$ which can be considered arbitrarily large and remembering that $C$ and $\epsilon$ have no dependence on $t$, we are allowed to take the limit infimum as $t\to+\infty$, to infer that 
			$$\liminf_{t \to -\infty} \frac{\lambda_{\D}(K, \phi_{-t}(K))}{t} \ge\frac{\epsilon}{\pi}.$$
			Once again, the choice of $\epsilon$ is completely random, so allowing $\epsilon\to+\infty$, we obtain the desired limit.
		\end{proof}
		
	\end{proposition} 
	
	\begin{remark}
		In this section, we almost exclusively examined the hyperbolic length and the extremal length along backward orbits resulting from non-elliptic semigroups. Following the methods introduced in Sections \ref{sec:Hyp_Length} and \ref{sec:Extremal_length}, respectively, we may extract analogous results in the case of an elliptic semigroup. We omit this study for the sake of not being redundant. 
	\end{remark}

\end{document}